\documentclass{article}
\usepackage{a4wide}
\usepackage{amsmath,amssymb}
\usepackage{arydshln}
\usepackage{algorithmic}
\usepackage{nicefrac}
\usepackage{graphicx,subcaption,booktabs}

\usepackage{tikz,pgfplots}
\usetikzlibrary{external}
\newtheorem{example}{Example}
\newtheorem{theorem}{Theorem}
\newtheorem{lemma}{Lemma}
\newenvironment{proof}{\textbf{Proof.\ }}{\hfill$\Box$~\\}

\title{Time integration of finite element models with nonlinear frequency dependencies}
\author{Elke Deckers \and Stijn Jonckheere \and Karl Meerbergen}
\date{\today}

\begin{document}
\maketitle

\begin{abstract}
The analysis of sound and vibrations is often performed in the frequency domain, implying the assumption of steady-state behaviour and time-harmonic excitation. External excitations, however, may be transient rather than time-harmonic, requiring time-domain analysis. Some material properties, e.g.\ often used to represent for damping treatments, are still described in the frequency domain, which complicates simulation in time.
In this paper, we present a method for the linearization of finite element models with nonlinear frequency dependencies.
The linearization relies on the rational approximation of the finite element matrices by the AAA method.
We introduce the Extended AAA method, which is classical AAA combined with a degree two polynomial term
to capture the second order behaviour of the models.
A filtering step is added for removing unstable poles.
\end{abstract}

\section{Introduction}

Classical vibro-acoustic analysis relies on a description of the model in the frequency domain. Passive damping and absorptive materials, such as visco-elastic \cite{Bagley1983}, porous \cite{Johnson1987,Champoux1991} and poro-elastic \cite{Biot1956a,Biot1956b} materials exhibit viscous and thermal damping mechanisms, which result in complex, frequency-dependent behaviour.  Many different descriptions to account for their complex behaviour can be found in literature, e.g.\ see \cite{Allard2009, Deckers2015}. Time domain analysis of such systems has gained quite some attention in the context of auralisation \cite{Vorlander2010,Jagla2012}, virtual sensing \cite{vandeWalle2018,vanophem2019} and inverse characterization \cite{Lewandowski2010}. Starting from the corresponding time-domain description of these materials is not always so straightforward as convolutions are required to account for the constitutive relationships. One solution is to make use of a recursive convolution which requires storage of only a minimum number of time steps per variable \cite{Semlyen1975,Kuether2016,Dragna2015}.
The goal of this paper is to present a method that allows to simulate nonlinear frequency dependent models directly in the time domain.

In its most general form, the frequency dependent model can be expressed as
\begin{eqnarray}
A(\omega) \hat{x} & = & \hat{b}(\omega), \label{eq:sys}
\end{eqnarray}
where $x$ is the state vector of dimension $n$ and $\omega$ is the angular frequency.
In classical FE analysis, the system matrix is quadratic in $\omega$. In this paper, we consider $A(\omega)$ that is not polynomial in $\omega$.
In fact, we assume that the system in the frequency domain can be written as the holomorphic decomposition or split form
\begin{eqnarray}
(A_0 + s A_1 + s^2 A_2 + A_{-1} g_1(s) + \cdots + A_{-m} g_m(s)) \hat{x} & = & \hat{b}(\omega) \label{eq:our-sys}
\end{eqnarray}
where $g_i$ is a scalar function, holomorphic on the imaginary axis, and the number of nonlinear terms,
$m$, is not large, i.e., a few dozen at most. 
%The norm $\|K\|_1 = 1.89 \cdot 10^9$ is much larger than $\|M\|_1 = 4.63 \cdot 10^{-4}$ and $\|C\|_1=5.96 \cdot 10^{-4}$.

For the simulation in the time domain of \eqref{eq:sys}, we aim to find a linear model
\begin{eqnarray}
  -\mathbf{E} \frac{d \mathbf{x}}{d t} + \mathbf{A} \mathbf{x} & = & \mathbf{b}(t)\quad,\quad t\geq 0, \\
           \mathbf{x}(0) & = & \mathbf{x}_0, \nonumber
\end{eqnarray}
such that the Laplace transform of $\mathbf{x}$ and $\hat{x}$ are `strongly' related.
First, we will derive a linear model in the frequency domain that approximates \eqref{eq:sys}.
The link with the time domain is then straightforward.
In order to form a linear model, we will use ideas from the solution of nonlinear eigenvalue problems.
For the relation with the time domain, we introduce the Laplace variable $s=i\omega$ and represent $A$ as a function of $s$.
In the sequel, we therefore use $s$ instead of $\omega$.
When $A(s)$ is a matrix polynomial or a rational matrix, i.e., the entries of $A$ are polynomials or
rational functions, there always are $\mathbf{E}$ and $\mathbf{A}$, $\mathbf{b}$ and $\mathbf{c}$ so that
\begin{eqnarray}
  -s \mathbf{E} \hat{\mathbf{x}} + \mathbf{A} \hat{\mathbf{x}} & = & \hat{\mathbf{b}} \label{eq:lin s}
\end{eqnarray}
for $s\in\imath\mathbb{R}$, and a way to extract $\hat{x}$ from $\hat{\mathbf{x}}$.
Model \eqref{eq:lin s} is called a \emph{linearization} of \eqref{eq:sys}.

When $A(s)$ is not a matrix polynomial or rational matrix, $A(s)$ is approximated by a rational matrix and the latter is then written in linear pencil form.
Each $g_j$ for $j=1,\ldots,m$ is approximated by a rational function.
In the literature, there are several techniques proposed for solving nonlinear eigenvalue problems.
In \cite{suba11}, a Pad\'e approximation is suggested.
Potential theory was used in NLEIGS \cite{gvmm14}.
The choice of poles and interpolation points is determined by the selection of the domain $\Sigma$ and a singularity set.
A rational approximation based on contour integration is proposed by \cite{saem19} and uses a basis of rational monomials.

In this work, we use the AAA method for rational approximation \cite{nast19} and the AAA-least squares variant \cite{cotr21}.
An efficient method for approximating $m$ functions was presented in \cite{FASTAAA}.
The set-valued AAA method uses the same ideas with the aim of developing a compact linearization of \eqref{eq:sys}  \cite{lmpv21}\cite{lime18} .
The weighted AAA method \cite{gnpt20} presents an alternative stopping criterion taking into account the relative contribution of each $g_j$ to $A(s)$.
The linearization from \cite{lmpv21} is complex valued and is not strong, in that it has an eigenvalue at infinity and the poles of the rational
function can be eigenvalues too.
These observation are important in the context of time integration.

For time stepping, it is important that $\mathbf{E}$, $\mathbf{A}$ and $\mathbf{b}$ are real.
Also, spurious eigenvalues arising from the linearization should be avoided.
In this paper, we use the notion of \emph{real function} for functions that are the result of the analytical extension of a real function
to the complex plane.
In this case, the rational approximation has real coefficients.
We denote the set of real functions by $\mathbb{F}_R$.
For any $g\in\mathbb{F}_R$, we have $g(\overline{z})=\overline{g}(z)$.

This paper is organised as follows. Section~\ref{sec:aaa} shows the AAA method for real functions.
The support points are selected as complex conjugate pairs. This needs some minor modifications to the
original algorithm from \cite{nast19} \cite{lmpv21}.
To ensure stability for time integration, we add filtering of unstable poles, using the AAA least squares method \cite{cotr21}.
We also introduce the \emph{extended AAA} algorithm, which builds a rational plus polynomial approximation.
The polynomial part reflects the quadratic dependence of finite element matrices on the frequency.
Section~\ref{sec:bary} proposes new real linearization pencils.
For AAA, we also present a linearization that does not have an infinite eigenvalue.
In Section~\ref{sec:time}, we discuss the connection of the nonlinear frequency model and the linear model in the time domain.
%Section~\ref{sec:mor} is devoted to a Greedy method for model order reduction  of large scale models.
Numerical examples are shown in Section~\ref{sec:numex}. Conclusions are given in Section~\ref{sec:conc}.

\section{AAA algorithm for real functions}\label{sec:aaa}

The AAA algorithm, pronounced as Triple A, is an iterative method that builds an approximation to a function $g$, expressed in barycentric rational Lagrange basis
\[
g(s) \approx r_d(s) = \frac{n_d(s)}{d_d(s)} = \frac{ \sum_{j=1}^d g(\sigma_j) \xi_j/(s-\sigma_j) }{\sum_{i=1}^d \xi_i/(s-\sigma_i)}.
%= [g(\sigma_1),\ldots,g(\sigma_d)]\cdot \Phi_d(s).
\]
The support points $\sigma_i$ and weights $\xi_i$ are selected iteratively with the aim to minimize an approximation error criterion.
The method works as follows.
Choose $n_Z$ test points $\{Z_j\}_{1\leq j\leq n_Z}$ in the region of interest.
Iteration $i$ consists of two steps.
In the first step, the support point $\sigma_i$ is chosen among the set $\{Z_j\}_{1\leq j\leq n_Z}\setminus\{\sigma_1,\ldots,\sigma_{i-1}\}$, so that
\[
\sigma_i = \text{argmax}_{s\in\{Z_j\}_{1\leq j\leq N}} |g(s) - n_d(s)/d_d(s)|.
\]
In the second step, the weights are chosen to minimize the linearized residual
$|g(s)d_d(s) - n_d(s) |$ for $s\in\{Z_j\}_{1\leq j\leq N}$.
The evaluation of the linearized residual in the sample points (minus the support points) leads to a linear least squares problem,
\[
  \left\| \begin{bmatrix}
g(Z_1) \sum_{j=1}^i \frac{\xi_i}{Z_1-\sigma_i}
 - \sum_{j=1}^i \frac{\xi_i g(\sigma_i)}{Z_1-\sigma_i}\\
    \vdots \\
g(Z_N) \sum_{j=1}^i \frac{\xi_i}{Z_1-\sigma_i}
- \sum_{j=1}^i \frac{\xi_i g(\sigma_i)}{Z_N-\sigma_i}.
\end{bmatrix}\right\|,
\]
which can be rewritten as
\[
  \left\| \begin{bmatrix}
  \frac{g(Z_1)-g(\sigma_1)}{Z_1-\sigma_1} & \cdots & \frac{g(Z_1)-g(\sigma_i)}{Z_1-\sigma_i} \\
    \vdots& & \vdots \\
  \frac{g(Z_{n_Z})-g(\sigma_1)}{Z_{n_Z}-\sigma_1} & \cdots & \frac{f(Z_{n_Z})-g(\sigma_i)}{Z_{n_Z}-\sigma_i}
\end{bmatrix} \begin{bmatrix} \xi_1 \\\vdots\\\xi_i \end{bmatrix}\right\|.
\]
Since the scaling of the weights does not alter the value of the rational approximation,
we can assume that the vector of weights is normalized.
As a result, the optimal weights correspond to the (normalized) right singular vector associated with the smallest singular value of
the Loewner matrix
\[
  L_i = \begin{bmatrix}
  \frac{g(Z_1)-g(\sigma_1)}{Z_1-\sigma_1} & \cdots & \frac{g(Z_1)-g(\sigma_i)}{Z_1-\sigma_i} \\
    \vdots & & \vdots\\
  \frac{g(Z_N)-g(\sigma_1)}{Z_N-\sigma_1} & \cdots & \frac{g(Z_N)-g(\sigma_i)}{Z_N-\sigma_i}
\end{bmatrix}.
\]
This procedure is iterated until the approximation error is sufficiently small for all points in $Z_N$.

\subsection{Real AAA}
The goal is to aim for a rational function that is a real function.
This means that $r(\overline{s})=\overline{r(s)}$.
This condition is satisfied if the support points as well as the weights are chosen as complex conjugate pairs.
In each iteration, we select a set of complex conjugate support points.
This means that when $\sigma_i$ is selected as a support point, also $\overline{\sigma}_i$ will be added to the set, i.e.,
if the support point is complex, two points are added in the same iteration.
The following describes the real AAA algorithm.
\begin{algorithmic}[1]
\STATE Given the set $Z_N$ and real function $g$.
\STATE $i=1$
\REPEAT
\STATE Find the next support point from
\[
\sigma_i = \text{argmax}_{s\in\{Z_j\}_{1\leq j\leq N}} |g(s) - r_d(s)|.
\]
\STATE Remove $\sigma_i$ and $\overline{\sigma}_i$ from the set $Z_N$.
\IF{$\mbox{Im}(\sigma_i)\neq0$}
\STATE Let $\sigma_{i+1}=\overline{\sigma}_i$.
\STATE $i=i+1$
\ENDIF
\STATE Set up the Loewner matrix $L_i$ and compute the smallest singular vector.
\STATE $i=i+1$
\UNTIL{the approximation error is small enough in all points in $Z_N$}
\end{algorithmic}

\if 0
For adding a support point at infinity,
the linearized residual, evaluated at $Z_\ell$, is written as
\begin{eqnarray*}
g(Z_\ell) \left( \xi_0 + \sum_{j=1}^i \frac{\xi_j}{Z_i-\sigma_j}\right) - \xi_0 g(\infty) - \sum_{j=1}^i \frac{\xi_j g(\sigma_j)}{Z_\ell-\sigma_j} & = & \\
(g(Z_\ell)-g(\infty)) \xi_0 + \sum_{j=1}^i \frac{g(Z_\ell)-g(\sigma_j)}{Z_\ell-\sigma_j}\xi_j.
\end{eqnarray*}
This leads to the following problem to be minimized:
\[
  \left\| \begin{bmatrix}
  (g(Z_1)-g(\infty)) & \frac{g(Z_1)-g(\sigma_1)}{Z_1-\sigma_1} & \cdots & \frac{g(Z_1)-g(\sigma_i)}{Z_1-\sigma_i} \\
    \vdots& \vdots & & \vdots \\
  (g(Z_N)-g(\infty)) & \frac{g(Z_N)-g(\sigma_1)}{Z_N-\sigma_1} & \cdots & \frac{g(Z_N)-g(\sigma_i)}{Z_N-\sigma_i}
\end{bmatrix} \begin{bmatrix} \xi_0 \\\vdots\\\xi_i \end{bmatrix}\right\|.
\]
\fi

The solution of the singular value problem will, in general, not generate weights with the desired complex conjugate structure of support points and weights, since the singular values also come in pairs
due to symmetry.
We therefore reformulate the problem as a real one as follows.
Reorder the support points $Z_1,\ldots,Z_N$ so that they come in complex conjugate pairs.
Define
\begin{eqnarray*}
T & = & \frac{1}{\sqrt{2}} \begin{bmatrix} 1 & 1 \\ \frac{1}{i} & -\frac{1}{i} \end{bmatrix}\quad,\quad T^{-1} = T^*, \\
S & = & \begin{bmatrix} 1 & i \\ 1 & -i \end{bmatrix},\\
S^{-1} & = & \frac{1}{2} \begin{bmatrix} 1 & 1 \\ \frac{1}{i} & -\frac{1}{i} \end{bmatrix}.
\end{eqnarray*}
Then define the unitary transformation $\mathbb{T}\in\mathbb{C}^{N\times N}$ such that there is $1$ on the
main diagonal when $Z_i$ is real, or $T$ when $Z_i$ is complex.
We also define $\mathbb{S}\in\mathbb{C}^{d\times d}$ as a block diagonal matrix with $1$ on the $i$th element of the main diagonal when $\sigma_i$ is real and $S$ otherwise.
The optimization problem is then written as
\[
\min \left\|
\begin{bmatrix}
\frac{g(Z_1) - g(\sigma_1)}{Z_1-\sigma_1} & \cdots & \frac{g(Z_1) - g(\sigma_d)}{Z_1-\sigma_d} \\
    \vdots & & \vdots \\
\frac{g(Z_{n_Z}) - g(\sigma_1)}{Z_{n_Z}-\sigma_1} & \cdots & \frac{g(Z_{n_Z}) - g(\sigma_d)}{Z_{n_Z}-\sigma_d}
\end{bmatrix}\mathbb{S} \cdot \mathbb{S}^{-1} \begin{pmatrix} \xi_1 \\\vdots \\\xi_d \end{pmatrix} \right\|
\]
Multiplying with $\mathbb{T}$ on the left does not change the norm:
\[
\min \left\| \mathbb{T}
\begin{bmatrix}
\frac{g(Z_1) - g(\sigma_1)}{Z_1-\sigma_1} & \cdots & \frac{g(Z_1) - g(\sigma_d)}{Z_1-\sigma_d} \\
    \vdots & & \vdots \\
\frac{g(Z_{n_Z}) - g(\sigma_1)}{Z_{n_Z}-\sigma_1} & \cdots & \frac{g(Z_{n_Z}) - g(\sigma_d)}{Z_{n_Z}-\sigma_d}
\end{bmatrix}\mathbb{S} \cdot \mathbb{S}^{-1} \begin{pmatrix} \xi_1 \\\vdots \\\xi_d \end{pmatrix} \right\|
\]
\begin{theorem}
The matrix
\[
\mathbf{L} = \mathbb{T}
\begin{bmatrix}
\frac{g(Z_1) - g(\sigma_1)}{Z_1-\sigma_1} & \cdots & \frac{g(Z_1) - g(\sigma_d)}{Z_1-\sigma_d} \\
    \vdots & & \vdots \\
\frac{g(Z_{n_Z}) - g(\sigma_1)}{Z_{n_Z}-\sigma_1} & \cdots & \frac{g(Z_{n_Z}) - g(\sigma_d)}{Z_{n_Z}-\sigma_d}
\end{bmatrix}\mathbb{S}
\]
is real.
\end{theorem}
\begin{proof}
We give a sketch of the proof.
Let $\sigma_i$ be real and $Z_j$ be complex and $Z_{j+1} = \overline{Z_j}$.
Then the column related to $\sigma_i$ is of the form
\[
\begin{matrix}
\cdots & \frac{G_j - g_i}{Z_j - \sigma_i} & \cdots \\
\cdots & \frac{\overline{G_j} - g_i}{\overline{Z_j} - \sigma_i} & \cdots
\end{matrix}
\]
Multiplication on the left with $\mathbb{T}$ produces
\[
\begin{matrix}
\cdots & \sqrt{2} \mbox{Re}\left( \frac{G_j - g_i}{Z_j - \sigma_i} \right) & \cdots \\
\cdots & \sqrt{2} \mbox{Im}\left( \frac{G_j - g_i}{Z_j - \sigma_i} \right) & \cdots
\end{matrix}
\]
This column is thus real.

Assume that $\sigma_i$ is complex and $\sigma_{i+1}=\overline{\sigma_i}$.
Similarly, assume that $Z_j$ is complex and $Z_{j+1} = \overline{Z_j}$.
First, multiply the Loewner matrix on the right with $\mathbb{S}$.
Then the two columns related to $\sigma_i$ and $\sigma_{i+1}$ are transformed into pairs of lines of the form
\begin{equation}\label{eq:proof-complex}
\begin{matrix}
\cdots & \frac{G_j - g_i}{Z_j - \sigma_i} + \frac{G_j - \overline{g_i}}{Z_j - \overline{\sigma_i}} &
 -\frac{1}{i}\left( \frac{G_j - g_i}{Z_j - \sigma_i} - \frac{G_j - \overline{g_i}}{Z_j - \overline{\sigma_i}}\right) & \cdots \\
\cdots & \frac{\overline{G_j} - g_i}{\overline{Z_j} - \sigma_i} + \frac{\overline{G_j} - \overline{g_i}}{\overline{Z_j} - \overline{\sigma_i}} &
 -\frac{1}{i}\left( \frac{\overline{G_j} - g_i}{\overline{Z_j} - \sigma_i} - \frac{\overline{G_j} - \overline{g_i}}{\overline{Z_j} - \overline{\sigma_i}}\right) & \cdots
\end{matrix}
\end{equation}
Now we multiply on the left with $\mathbb{T}$.
This transforms the two lines from above to
\[
\begin{matrix}
\cdots & \frac{\sqrt{2}}{2}\mbox{Re}\left( \frac{G_j - g_i}{Z_j - \sigma_i} + \frac{G_j - \overline{g_i}}{Z_j-\overline{\sigma_i}} \right) & 
                        \frac{\sqrt{2}}{2}\mbox{Im}\left( \frac{G_j - \overline g_i}{Z_j - \overline \sigma_i} - \frac{G_j - {g_i}}{Z_j-{\sigma_i}} \right) \\
\cdots & \frac{\sqrt{2}}{2}\mbox{Im}\left( \frac{G_j - g_i}{Z_j - \sigma_i} + \frac{G_j - \overline{g_i}}{Z_j-\overline{\sigma_i}} \right) &
                        -\frac{\sqrt{2}}{2}\mbox{Re}\left( \frac{G_j - \overline g_i}{Z_j - \overline \sigma_i} - \frac{G_j-{f_i}}{Z_j-{\sigma_i}} \right)
\end{matrix}
\]

If $\sigma_i$ is complex, but $Z_j$ is real, then the first line in Eq.~\eqref{eq:proof-complex} is real:
\[
\begin{matrix}
\cdots & 2 \mbox{Re}\left( \frac{G_j - g_i}{Z_j - \sigma_i} \right) &
 - 2\mbox{Im}\left( \frac{G_j - g_i}{Z_j - \sigma_i} \right) & \cdots \\
\end{matrix}
\]
\end{proof}
Let the smallest singular vector of $\mathbf{L}$ be $w$, then the weights are obtained as $\mathbb{S} w$.
Since $\mathbf{L}$ is real, $w$ is real, and therefore, $\mathbb{S} w$ has the desired complex conjugate structure.

\subsection{Filtered AAA}\label{sec:AAA-LS}
There is, in general, no guarantee that the AAA method leads to stable poles.
Our numerical experiments show that often all poles are stable, but a few unstable poles of large modulus sometimes appear.
These poles of large modulus are usually not so important and can therefore be removed.

We use the AAA-least squares method from \cite{cotr21} for this purpose.
In the remainder of the text, we call this approach \emph{Filtered AAA}, or F-AAA, for short.
The barycentric form does not allow to remove poles explicitly. Therefore, another rational basis is used.

Let the poles obtained by AAA be $\mu_1,\ldots,\mu_d$ and let only the poles $\mu_1,\ldots,\mu_k$ with $k\leq d$ be stable.
Then, we build another rational approximation with the stable poles only, using weighted partial fractions.
Let
\[
  \gamma_0 + \gamma_1 \frac{\xi_1}{s-\mu_1} + \cdots + \gamma_k\frac{\xi_k}{s-\mu_k}
\]
be the weighted partial fractions approximation to $g(s)$.
We assume that the poles are simple.
This makes sense, since,
it will be clear from the proof of Theorem~\ref{th:eig} in \S\ref{sec:bary} that the poles are simple, because the dual linear basis of the rational functions has full rank.
The values $\xi_1,\ldots,\xi_k$ are strictly positive real weights, chosen so that $\xi_i/|s-\mu_i|$ is not larger than one in the sample points,
and the values $\gamma_0,\ldots,\gamma_k$ are the coefficients related to $g(s)$.
The coefficients $\gamma_0,\ldots,\gamma_k$ are determined from the linear least squares problem
\[
\min \sum_{s\in\mathbf{Z}} \left|\gamma_0 + \gamma_1\frac{\xi_1}{s-\mu_1} + \cdots + \gamma_k\frac{\xi_k}{s-\mu_k} - g(s)\right|^2.
\]
In matrix form, the overdetermined linear system becomes
\[
  \begin{bmatrix} 1 & \frac{\xi_1}{Z_1-\mu_1} & \cdots & \frac{\xi_k}{Z_1-\mu_k} \\ \vdots & & & \vdots \\
    1 & \frac{\xi_1}{Z_N-\mu_1} & \cdots & \frac{\xi_k}{Z_N-\mu_k} \end{bmatrix} \begin{bmatrix} \gamma_1 \\ \vdots \\\gamma_k \end{bmatrix} =
    \begin{bmatrix} g(Z_1) \\ \vdots \\ g(Z_N) \end{bmatrix}.
\]
By the scaling, the columns of the matrix of the overdetermined system, have infinity norm at most one.
By grouping complex conjugate poles and complex conjugate sample points together, a transformation maps the matrix to a real one and a real right-hand side, so that, after
backtransformation, the coefficients $\gamma_0,\ldots,\gamma_k$ are such that the rational approximation to $g$ is a real function.

In the case of a situation where some poles are clustered, the proposed basis might be ill-conditioned and it may be better to use the following
basis, which is the inverse of Newton polynomials:
\[
  \gamma_0 + \gamma_1\frac{\xi_1}{s-\mu_1} + \gamma_2\frac{\xi_2}{(s-\mu_1)(s-\mu_2)} \cdots + \gamma_k\frac{\xi_k}{(s-\mu_1)\cdots(s-\mu_k)}.
\]
We order the poles so that $\mu_1,\ldots,\mu_{k_r}$ are real and the poles
\[
  \mu_{k_r+1},\overline{\mu}_{k_r+1}\ldots,\mu_{k_r+k_c},\overline{\mu}_{k_r+k_c}
\]
are complex, where $k=k_r+2k_c$.
Then, define the basis functions
\begin{eqnarray*}
\phi_0 & = & 1, \\
\phi_{j} & = & \frac{\xi_{j}}{(s-\mu_1)\cdots(s-\mu_{j})} , \\
\tilde{\phi}_j & = & \frac{\xi_{j}}{(s-\overline{\mu_1})\cdots(s-\overline{\mu_{j}})}.
\end{eqnarray*}
The rational approximation is written as
\begin{multline*}
\gamma_0 \phi_0(s) + \gamma_1 \phi_1(s) + \cdots + \gamma_{k_r} \phi_{k_r}(s) \\
    + \gamma_{k_r+1} \phi_{k_r+1}(s) + \cdots + \gamma_{k_r+k_c} \phi_{k_r+k_c}(s) \\
    + \overline{\gamma_{k_r+1}} \tilde{\phi}_{k_r+1}(s) + \cdots + \overline{\gamma_{k_r+k_c}} \tilde{\phi}_{k_r+k_c}(s) .
\end{multline*}
For example, with poles $1$, and $2\pm3\imath$, we have the rational approximation
\[
\gamma_0 + \gamma_1\frac{\xi_1}{s-1} + \gamma_2\frac{\xi_2}{(s-1)(s-2-3\imath)} + \overline{\gamma_2}\frac{{\xi_2}}{(s-1)(s-2+3\imath)}.
\]
As for the partial fraction representation, the least squares problem can be transformed to a real valued least squares problem, which then leads to
the coefficients with complex conjugate structure
\[
  \gamma_0\ \gamma_1\ \cdots \gamma_{k_r}\ \gamma_{k_r+1} \cdots \gamma_{k_r+k_c}\ \overline{\gamma}_{k_r+1} \cdots \overline{\gamma}_{k_r+k_c}.
\]

\subsection{Set-valued AAA and weighted AAA}
For the approximation of \eqref{eq:our-sys} by a rational function, rational approximations to each $g_1,\ldots,g_m$ are developed.
It was observed in \cite{lmpv21} that this may lead to a large size linearization.
Therefore, the set-valued method was proposed, which builds $m$ approximations $r_1,\ldots,r_m$ to $g_1,\ldots,g_m$ simultaneously, where all $r_j$ use the
same support points and weights.
An alternative procedure was proposed in \cite{FASTAAA}.

The support points are chosen from a maximization criterion, i.e., the worst point among all test points is taken.
The weights, and therefore, the poles are chosen from the Loewner matrix in the Set-valued AAA method.
The Loewner matrix is obtained by vertically stacking the individual Loewner matrices of the $m$ functions.
The scaling of functions may therefore have an impact on the weights.
We  do not discuss the technical details in this paper, but refer to \cite{FASTAAA,lmpv21,gnpt20}.

The Set-valued AAA method builds Loewner matrices for $g_1/\|g_1\|_Z,\ldots,g_m/\|g_m\|_Z$ and uses the stop criterion
\begin{equation}\label{eq:stop aaa max}
  \|r_j(s) - g_j(s)\|_Z \leq \tau \|g_j(s)\|_Z,
  \quad\text{for}\quad j=1,\ldots,m,
\end{equation}
where $\tau$ is a  selected tolerance.
The criterion in \eqref{eq:stop aaa max} is also used to select the next support point by selecting $Z_i$ that maximizes
\[
  |r_j(Z_i) - g_j(Z_i)| / \|g_j(s)\|_Z,
  \quad\text{for}\quad j=1,\ldots,m.
\]
The disadvantage of this choice is that all functions are treated in the same way, even if they do not contribute much to $A(s)$, i.e., the selection of the
support points and the weight is independent of the contribution to $A(s)$.

Another proposal, Weighted AAA, was suggested in \cite{gnpt20}, that presents an adapted way to select the support points.
We use Weighted AAA in the numerical experiments.
The functions are now scaled so that
\[
  |r_j(Z_i) - g_j(Z_i) \|A_{-j}\|_2| ,
  \quad\text{for}\quad j=1,\ldots,m,
\]
is minimized by the same greedy optimization procedure as set-valued AAA.
A tuned stop criterion takes into account the contribution of each $g_j(s)$ in $A(s)$.
Indeed, if $\|A_{-j}\|$ is small, the accuracy of the rational approximation of $g_j$ should not be high.
The main idea is to find $R(s)$ that approximates $A(s)$ with a relative error
\[
  \|A(s) - R(s)\|_Z \leq \tau \|A(s)\|_Z.
\]
The difficulty is a cheap estimate of the norms.
A lower bound on $\|A\|_2$ can be found by a single matrix vector product with a randomly chosen vector $v$, i.e., a vector whose elements are drawn from a normal distribution on $[-1,1]$, e,g.
and then using the estimate
\[
  \|A\|_F \approx \sqrt{n}\frac{\|A v\|_2}{\|v\|_2},
\]
see \cite{gukl95}.
We express $A_j v = Q w_j$ with $Q\in\mathbb{C}^{n\times 3+m}$ and $Q_j^*Q_j=I$.
Then, for any $s\in\mathbb{C}$,
\begin{eqnarray*}
\|A(s) v\|_2 & = & \|w_0 + s w_1 +s^2 w_2 + g_1(s) w_{-1} + \cdots + g_m w_{-m}\|_2. \\
\|A(s) v - R(s) v \|_2 & = & \|(g_1(s)-r_1(s)) w_{-1} + \cdots + (r_m(s)-g_m-r_m(s)) w_{-m}\|_2.
\end{eqnarray*}

\section{Linearizations of real rational functions}\label{sec:bary}
In this section, we represent a rational function, expressed in barycentric form or partial fraction form as the Schur complement of a linear pencil
\[
  L(s) = \begin{bmatrix} \alpha & a^T \\ c_0+sc_1 & C_0+sC_1 \end{bmatrix}
\]
with $\alpha\in\mathbb{C}$, $a,c_0,c_1\in\mathbb{C}^{\tilde{d}}$ and $C_0,C_1\in\mathbb{C}^{\tilde{d}\times \tilde{d}}$.
The Schur complement is
\begin{equation}\label{eq:schur complement}
  r(s) = \alpha - a^T (C_0+sC_1)^{-1}(c_0+sc_1).
\end{equation}
The factors
\[
  \Phi = -(C_0+sC_1)^{-1}(c_0+sc_1)
\]
are a set of rational basis functions and the vector $a$ contains the coefficients, i.e., $r = \alpha+ a^T \Phi$.
It satisfies
\begin{equation}\label{eq:rat from schur complement}
  \begin{bmatrix} \alpha & a^T \\ c_0+sc_1 & C_0+sC_1 \end{bmatrix} \begin{pmatrix} 1 \\ \Phi(s) \end{pmatrix} = \begin{pmatrix} r(s) \\ 0 \end{pmatrix}.
\end{equation}
The linearization $L(s)$ can now be used to linearize a rational matrix expressed in the basis $\Phi$.

We further assume that $\begin{bmatrix} c_0+sc_1 & C_0+sC_1 \end{bmatrix}$ is a full rank matrix for all $s\in\mathbb{C}$.
We also assume that $C_1$ has full rank. This implies that the rational basis functions $\Phi$ do not have a pole at infinity, such that the constant function $1$ cannot be spanned by the elements of
$\Phi$.
Note that the original AAA linearization has a singular $C_1$.
Indeed, the sum of the rational basis functions leads to the constant function $1$.
We will further show a linearization that removes the pole at infinity.
We can therefore assume that $C_1$ is nonsingular.

\begin{lemma}\label{le:schur transformation}
Let a rational function be defined by the Schur complement of $L(s)$.
Let $Y\in\mathbb{C}^{\tilde{d}\times \tilde{d}}$, and $Z\in\mathbb{C}^{\tilde{d}\times \tilde{d}}$
be full rank matrices.
Then, the Schur complement of $L(s)$ is the same as the Schur complement of
\[
\begin{bmatrix}
\tilde{\alpha} & \tilde{a}^T \\ \tilde{c}_0+s\tilde{c}_1 & \tilde{C}_0+s\tilde{C}_1 \end{bmatrix} = \begin{bmatrix} 1 & 0 \\ 0 & Y\end{bmatrix}
\begin{bmatrix}
\alpha & a^T \\ c_0+sc_1 & C_0+sC_1 \end{bmatrix} \begin{bmatrix} 1 & 0 \\ z & Z\end{bmatrix}
\]
\end{lemma}
\begin{proof}
Multiplying \eqref{eq:rat from schur complement} on the left with $\begin{bmatrix} 1 & 0 \\ 0 & Y\end{bmatrix}$ does not change the right-hand side in \eqref{eq:rat from schur complement}.
Let
\[
\mathbf{\Phi}=\begin{pmatrix} 1 \\ \Phi\end{pmatrix}.
\]
The multiplication with $Z$, changes the basis to
\[
\tilde{\mathbf{\Phi}}= \begin{pmatrix} 1 \\ Z^{-1}(\Phi - z)\end{pmatrix}
\]
The rational function is
\[
[\alpha\ a^T]\mathbf{\Phi} = [\tilde\alpha\ \tilde{a}^T] \tilde{\mathbf{\Phi}},
\]
where $[\tilde\alpha\ \tilde{a}^T]$ are the coefficients for the new basis $\tilde{\mathbf{\Phi}}$.
\end{proof}

\subsection{Linearization of rational matrices}
Consider the rational matrix
\begin{equation}\label{eq:rat approx}
  R(s) = A_0+sA_1+s^2A_2 + \sum_{j=1}^m A_{-j} r_j(s) \in\mathbb{C}^{n\times n},%(\alpha_j - a_j^T (C_0+sC_1)^{-1} (c_0+sc_1))
\end{equation}
where the rational function $r_j$ is expressed as
\[
  r_j(s) = \alpha_j + a_j^T \Phi(s),
\]
i.e., the basis $\Phi$ is the same for all rational functions.
\begin{theorem}
The system
\begin{eqnarray*}
R(s) \hat{x} & = &  \hat{b} \\
\end{eqnarray*}
with $R(s)$ given by \eqref{eq:rat approx},
is related to the linear system
\begin{eqnarray}\label{eq:rat lin}
  (\mathbf{A} - s \mathbf{E}) \hat{\mathbf{x}} & = & e_1 \otimes \hat{b},
\end{eqnarray}
with
\begin{eqnarray*}
\mathbf{L}(s) = \mathbf{A} - s \mathbf{E} & = & \left[\begin{array}{cc|c} \tilde{A}_0+s A_1 & s A_2 & a_1^T \otimes A_{-1} + \cdots + a_m^T \otimes A_{-m} \\
                           s I & -I & 0 \\\hline
                           (c_0+sc_1)\otimes I &0 & (C_0+sC_1)\otimes I
                \end{array}\right],\\
\tilde{A}_0 & = & A_0 + \alpha_1 A_{-1} + \cdots + \alpha_m A_{-m} \\
\hat{\mathbf{x}} & = & \mathbf{\Phi}(s) \otimes \hat{x}(s),\\
\mathbf{\Phi}(s) & = & \begin{pmatrix} 1 \\ s \\ \Phi(s) \end{pmatrix}.
\end{eqnarray*}
\end{theorem}
\begin{proof}
Working out the first row of $(\mathbf{A} - s \mathbf{E}) \hat{\mathbf{x}} = e_1 \otimes \hat{b}$ immediately leads to $R(s) \hat{x} = \hat{b}$.
The other rows impose the structure of $\hat{\mathbf{x}}$.
The second block row corresponds to the linearization of the quadratic term.
We also use $(C_0+C_1)\Phi+c_0+sc_1=0$.
\end{proof}
We now make a statement about the eigenvalues of the pencil $\mathbf{L}(s)$.

Let $R(s)$ be a rational matrix with the set of poles $\Xi$.
An eigenvalue of $R(s)$ is a value $s\in(\mathbb{C}\cup\{\infty\})\setminus \Xi$ for which $\det(R(s))=0$.
\begin{theorem}\label{th:eig}
The eigenvalues of $R(s)$ are eigenvalues of $\mathbf{L}(s)$, i.e., if $\det(R(s))=0$ then it follows that $\det(\mathbf{L}(s))=0$.
Let $[c_0+sc_1\ C_0+sC_1]$ have full rank for all $s\in\mathbb{C}\cup\{\infty\}$.
Then, the eigenvalues of $\mathbf{L}(s)$, that are not poles of $R(s)$, are eigenvalues of $R(s)$.

In addition, a pole $s$ of $R(s)$ is an eigenvalue of $\mathbf{L}$ iff
\[
\begin{bmatrix} a_1^T \otimes A_{-1} + \cdots + a_m^T \otimes A_{-m} \\  (C_0+sC_1)\otimes I
\end{bmatrix}
\]
does not have full column rank.
\end{theorem}
\begin{proof}
First, note that
\[
  M-sN = \begin{bmatrix} s & -1 & 0 \\ c_0+s c_1 & 0 & C_0+sC_1 \end{bmatrix}
\]
has full rank iff $[c_0+sc_1\ C_0+sC_1]$ has full rank.
Let $R(\lambda)x=0$ with $x\neq 0$, then
\[
  (\mathbf{A}-\lambda \mathbf{E}) \mathbf{\Phi}(\lambda) \otimes x =0.
\]
Now, to prove that any eigenvalue of $\mathbf{A}-s\mathbf{E}$ is an eigenvalue of $R(s)$, note that, since $M-sN$ has full rank,
null vectors of $M-sN$ always are a multiple of $\mathbf{\Phi}(s)$.
The eigenvectors of $\mathbf{A}-s\mathbf{E}$ must therefore have
the form $\mathbf{\Phi}(\lambda) \otimes x$.
It follows that then $R(\lambda)x=0$, for $\lambda$ for which $\|R(\lambda)\|$ is finite.

Since $[c_0+sc_1\ C_0+sC_1]$ has full rank for all $s\in\mathbb{C}\cup\{\infty\}$, the eigenvalues of $C_0+sC_1$ must be simple.
Let $\lambda$ be an eigenvalue of $C_0+sC_1$ and, therefore, be a pole of $R(s)$ with associated eigenvector $z$.
Since $M-\lambda N$ is full rank, multiples of
\[
 \mathbf{z}=\begin{pmatrix} 0 \\ 0 \\ z\end{pmatrix}
\]
are the only nullvectors of $M-\lambda N$.
Therefore, $(\mathbf{A}-\lambda \mathbf{E}) \mathbf{x}=0$ is only possible iff $\mathbf{x} = \mathbf{z}\otimes x$,
and
\[
\begin{bmatrix} a_1^T \otimes A_{-1} + \cdots + a_m^T \otimes A_{-m} \\  (C_0+\lambda C_1)\otimes I
\end{bmatrix} (z\otimes x) = 0.
\]
\end{proof}
Matrix $M-sN$ is called a dual linear basis of $\mathbf{\Phi}(s)$, i.e., $(M-sN)\mathbf{\Phi}=0$.
From this theorem, it follows that it is important that $[c_0+sc_1\ C_0+sC_1]$ or $M-sN$ has full rank for $s\in\mathbb{C}\cup\{\infty\}$.
We will therefore aim for linearizations $L$ with this property.
A dual linear basis with this property is called a full rank dual basis.

\subsection{Linearization for barycentric form}
Given $d$ support points $\sigma_j,j=1,\ldots,d$ and associated weights $\xi_j$, then, the barycentric Lagrange basis is
\[
\phi_j = \frac{\xi_j/(s-\sigma_j)}{\sum_{i=1}^d \xi_i/(s-\sigma_i)} .
\]
%We define $\gamma = \sum_{i=1}^d \xi_i/(s-\sigma_i)$.
Due to the scaling factor in the denominator of $\phi_j$, we have that $\sum_{j=1}^d \phi_j \equiv 1$.
It makes sense to assume that all $\xi_1,\ldots,\xi_d$ are nonzero otherwise some support points can be removed.
A rational function $r(s)$ is expressed as
\[
r(s) = r(\sigma_1) \phi_1(s) + \cdots + r(\sigma_d) \phi_d(s) .
\]

Now define $\Phi = [\phi_1,\ldots,\phi_d]^T$, and
\begin{equation}\label{eq:rat-schur}
{L}_1 = \begin{bmatrix}
0 & a^T \\ c_0+sc_1 & C_0+sC_1 \end{bmatrix} = \begin{bmatrix}
0 & r(\sigma_1) & \cdots & r(\sigma_d) \\
-1 & 1 & \cdots & 1 \\
0 & \xi_2(s-\sigma_1) & \xi_1(\sigma_2-s) & \\
\vdots &         & \ddots & \ddots & \\
0 &  & \xi_d(s-\sigma_{d-1}) & \xi_{d-1}(\sigma_d-s)
\end{bmatrix}
\end{equation}
then, we have that $(C_0+sC_1)\Phi = e_1$.
The first row imposes the scaling so that $\phi_1+\cdots+\phi_d=1$.
The other rows lead to the ratio of two successive basis functions.

Matrix \eqref{eq:rat-schur} is called a linearization of $r(z)$.
One difficulty with the formulation \eqref{eq:rat-schur} is that
$[c_1\ C_1]$ is rank deficient. In fact, $c_1=0$ and $C_1$ is singular. This introduces a linearization with an infinite eigenvalue.
For this reason, $\begin{bmatrix} c_0+sc_1 & C_0+sC_1 \end{bmatrix}$ is not a full rank dual linear basis of $\Phi$.
This may have some disadvantages in the numerical methods, as we will discuss later.

We now derive an alternative linearization by applying Lemma~\ref{le:schur transformation} with the transformations
\[
Y=I \quad\text{and}\quad Z = \begin{bmatrix} 1 & 0 \\ 0 & 1 & 0 \\ & & \ddots \\ 1 & -1 & \cdots & -1 & 1 \end{bmatrix}, z=0.
\]
The basis after transformation is $Z^{-1}{\Phi} = [1,\phi_1,\ldots,\phi_{d-1},0]$.
The linearization becomes
{\small
\[
\begin{bmatrix}
r(\sigma_d) & r(\sigma_1)-r(\sigma_d) & \cdots & \cdots & r(\sigma_{d-1})-r(\sigma_d) & r(\sigma_d) \\
0 & 0 & \cdots & \cdots & 0 & 1 \\
0 & \xi_2(s-\sigma_1) & \xi_1(\sigma_2-s) & \\
\vdots &         & \ddots & \ddots & \\
0 &  & \xi_{d-1}(s-\sigma_{d-2}) & \xi_{d-2}(\sigma_{d-1}-s) \\
\xi_{d-1}(\sigma_d-s) &  \xi_{d-1}(s-\sigma_d) & \cdots & \xi_{d-1}(s-\sigma_d) & \xi_{d-1}(s-\sigma_d) + \xi_d(s-\sigma_{d-1}) &  \xi_{d-1}(s-\sigma_d)
\end{bmatrix}
\]
}
The last column and second row do not contribute and can be removed:
\begin{equation}\label{eq:rat-schur-strong}
\begin{bmatrix}
r(\sigma_d) & r(\sigma_1)-r(\sigma_d) & \cdots & \cdots & r(\sigma_{d-1})-r(\sigma_d) \\
0 & \xi_2(s-\sigma_1) & \xi_1(\sigma_2-s) \\
\vdots &         & \ddots & \ddots & \\
0 &  & & \xi_{d-1}(s-\sigma_{d-2}) & \xi_{d-2}(\sigma_{d-1}-s) \\
\xi_{d-1}(\sigma_d-s) &  \xi_{d-1}(s-\sigma_d) & \cdots & \xi_{d-1}(s-\sigma_d) & \xi_{d-1}(s-\sigma_d) + \xi_d(s-\sigma_{d-1})
\end{bmatrix}
\end{equation}
This is a linearization for the basis $[1,\phi_1,\ldots,\phi_{d-1}]$, where $\phi_d=1-\phi_1-\cdots-\phi_{d-1}$.
\begin{lemma}\label{le:rat-schur-strong}
The last $d-1$ rows of \eqref{eq:rat-schur-strong} form a full rank dual linear basis to $[1,\phi_1,\ldots,\phi_{d-1}]^T$.
\end{lemma}
\begin{proof}
We have to prove that $[c_0+sc_1\ C_0 +s C_1]$ has full rank for all $s\in\mathbb{C}\cup\{\infty\}$.
Firstly, if $s=\sigma_d$, we have
\[
  \begin{bmatrix}
  0 & \xi_2(\sigma_d-\sigma_1) & \xi_1(\sigma_2-\sigma_d) \\
\vdots &         & \ddots & \ddots & \\
0 &  & & \xi_{d-1}(\sigma_d-\sigma_{d-2}) & \xi_{d-2}(\sigma_{d-1}-\sigma_d) \\
0 &  0 & \cdots & 0 & \xi_d(\sigma_d-\sigma_{d-1})
  \end{bmatrix}
\]
which has full rank, since the support points $\sigma_i$ are simple.
If $s\neq\sigma_d$, then the matrix has full rank iff the submatrix
\[
\begin{bmatrix}
\xi_2(s-\sigma_1) & \xi_1(\sigma_2-s) \\
 &   \ddots     & \ddots & \\
 & \xi_{d-1}(s-\sigma_{d-2}) & \xi_{d-2}(\sigma_{d-1}-s)
\end{bmatrix}
\]
has full rank.
It is easy to see that this is the case for any value of $s\in\mathbb{C}\cup\{\infty\}$.
This proves the lemma.
\end{proof}

\subsection{Real valued linearization for the barycentric form}
For real rational functions, we can derive linearization pencils with real matrices.
Recall the definition of a real function $r(s)$ as a function that satisfies the property $r(\overline{s})=\overline{r(s)}$.
For rational functions in barycentric form, we will ensure that the support points come in complex conjugate pairs.
As a result, it is not hard to see that the weights also form complex conjugate pairs.

We will use Lemma~\ref{le:schur transformation} to map the linearization \eqref{eq:rat-schur-strong} to a pencil with real valued coefficient matrices.
To derive this linearization, we use an alternative linearization to \eqref{eq:rat-schur}.
When there are only complex conjugate support points, then, we use the linearization
\[
{L}_2 = \begin{bmatrix}
  0 & \mathbf{r}(\sigma_1) & \cdots & \mathbf{r}(\sigma_{d-1})  \\
  -1 & \mathbf{e}^T & \cdots & \mathbf{e}^T \\
   0 & -L_{3,1}(s) & L_{1,3} \\
   \vdots & \ddots & \ddots \\
   0 && -L_{d-1,d-3}(s) & L_{d-3,d-1}  \\
   0 &&& -h_{d-1}(s)^T
  \end{bmatrix}
\]
with
\begin{eqnarray*}
\mathbf{r}(\sigma) & = & [r(\sigma), r(\overline{\sigma})], \\
\mathbf{e}^T & = & [1,1], \\
L_{k,l}(s) & = & \begin{bmatrix} \xi_k (\sigma_{l}-s) & 0 \\ 0 & \overline{\xi_k} (\overline{\sigma_l}-s) \end{bmatrix},\\
h_{l}(s)^T & = & \begin{bmatrix} \overline{\xi_l} (\sigma_{l}-s) & \xi_l (s - \overline{\sigma_l}) \end{bmatrix}.
\end{eqnarray*}
The linearization makes a coupling between pairs of basis functions, and in the last line a coupling is made between the last support point and its conjugate.
It is easy to see that
\[
  {L}_2(s) \Phi = e_1 r(s),
\]
which makes the relation with the rational function $r(s)$ and imposes the structure of the basis $\Phi$.
Note that the linearization does not have a full rank dual linear basis for $\Phi$ for the same reasons as \eqref{eq:rat-schur}.

Now, we transform the complex valued linearization pencil ${L}_2$ to a real one.
Define the following transformation matrices
\begin{equation}\label{eq:T}
T = \begin{bmatrix}
1 & 1 \\
\imath & -\imath
\end{bmatrix}\quad\text{and}\quad T^{-1} = \frac{1}{2} \begin{bmatrix} 1 & -\imath \\ 1 & \imath \end{bmatrix},
\end{equation}
and
\[
\mathbb{T}_R = \begin{bmatrix} 1 &  & & \\
      & T &  \\
      & & \ddots & \\
      & & & T \end{bmatrix}\quad\mbox{and}\quad
\mathbb{T}_L = \begin{bmatrix} 1 & \\ & 1 \\
      && T &  \\
      && & \ddots & \\
      && & & T \\ &&&&& -\imath\end{bmatrix},
\]
where $T$ is repeated as many times as there are complex conjugate pairs of support points.
We define the new basis
\[
{\Psi} = \mathbb{T}_R {\Phi}
\]
By construction, the entries of $\Psi$ are
\[
\Psi = \begin{bmatrix} 1 \\
    \left(\frac{\xi_{1}}{s-\sigma_{1}} + \frac{\overline{\xi_{1}}}{s-\overline{\sigma_{\ell+1}}}\right) \big/ \sum_{i=1}^d \frac{\xi_{i}}{s-\sigma_{i}}  \\
    \imath\left(\frac{\xi_{1}}{s-\sigma_{1}} - \frac{\overline{\xi_{1}}}{s-\overline{\sigma_{\ell+1}}}\right)  \big/ \sum_{i=1}^d \frac{\xi_{i}}{s-\sigma_{i}} \\
    \vdots
\end{bmatrix}
\]
Note that all elements of $\Psi$ lie in $\mathbb{F}_R$.
We define the linearization
\[
  {L}_3 = \mathbb{T}_L {L}_2 \mathbb{T}_R^{-1}.
\]
Following Lemma~\ref{le:schur transformation}, the basis has changed from ${\Phi}$ to ${\Psi}$,
and $L_2$ and $L_3$ have the same Schur complements.
What remains to show is that the coefficient matrices of $L_3$ are real.
We have a look at the different blocks of ${L}_3$.
\[
  {L}_3 = \begin{bmatrix}
0& \mathbf{r}(\sigma_{1}) T^{-1} & \cdots & \mathbf{r}(\sigma_{d-1}) T^{-1} \\
  -1  & \mathbf{e}^T T^{-1} & \cdots & \mathbf{e}^T T^{-1} \\
   & -T L_{3,1} T^{-1} & T L_{1,3} T^{-1} \\
           &         & & \ddots & \ddots
  0 & && -\imath h_{d-1}^T T^{-1} & \\
  \end{bmatrix},
\]
with
\begin{eqnarray*}
\tilde{\mathbf{e}}^T = \mathbf{e}^T T^{-1}  & = & [1\ 0], \\
\tilde{\mathbf{r}}(\sigma) = \mathbf{r}(\sigma) T^{-1}  & = & [\mbox{Re}(r(\sigma))\ \mbox{Im}(r(\sigma))], \\
\tilde{L}_{k,l} = T L_{k,l}(s) T^{-1} & = & \begin{bmatrix} \mbox{Re}(\xi_k \sigma_{l}) - \mbox{Re}(\xi_k) s & \mbox{Im}(\xi_k \sigma_{l}) + \mbox{Im}(\xi_k) s \\ -\mbox{Im}(\xi_k \sigma_{l}) + \mbox{Im}(\xi_k) s & \mbox{Re}(\xi_k \sigma_{l}) - \mbox{Re}(\xi_k) s \end{bmatrix},\\
\tilde{h}_l = -\imath h_{l}^T T^{-1} & = & \begin{bmatrix} \mbox{Im}(\xi_l \overline{\sigma_{l}}) - \mbox{Im}(\xi_l) s & \mbox{Re}(\xi_l \overline{\sigma_{l}}) - \mbox{Re}(\xi_l) s \end{bmatrix}.
\end{eqnarray*}

\begin{example}\label{ex:d=4}
For $d=4$ with complex support points $\sigma_1,\sigma_2=\overline{\sigma_1}\in\mathbb{R}$ and $\sigma_3,\sigma_4=\overline{\sigma_3}\in\mathbb{C}$,
\[
{L}_2  = \begin{bmatrix}
0 & r(\sigma_1) & \overline{r(\sigma_1)} & r(\sigma_3) & \overline{r(\sigma_3)} \\
-1 & 1 & 1 & 1 & 1 \\
0 & \xi_3(s-\sigma_1) & & \xi_1(\sigma_3-s) & 0 \\
0 & 0 &  \overline{\xi_3}(s-\overline{\sigma_1}) & 0 & \overline{\xi_1}(\overline{\sigma_3}-s) \\
0 & & & \overline{\xi_3}(s-\sigma_3) & \xi_3(\overline{\sigma_3}-s)
\end{bmatrix}
\]
and
\[
{L}_3 = \begin{bmatrix}
0 & \mbox{Re}(r(\sigma_1)) & \mbox{Im}(r(\sigma_1)) & \mbox{Re}(r(\sigma_3)) & \mbox{Im}(r(\sigma_3)) \\
-1 & 1 & 0 & 1 & 0 \\
0 & \mbox{Re}(\xi_3)s - \mbox{Re}(\xi_3\sigma_1) & \mbox{Im}(\xi_3) s -\mbox{Im}(\xi_3\sigma_1) & \mbox{Re}(\xi_1\sigma_3)-\mbox{Re}(\xi_1)s & \mbox{Im}(\xi_1\sigma_3)+\mbox{Im}(\xi_1)s \\
0 &  -\mbox{Im}(\xi_3)s + \mbox{Im}(\xi_3\sigma_1) & \mbox{Re}(\xi_3)s - \mbox{Re}(\xi_3\sigma_1) & -\mbox{Im}(\xi_1\sigma_3)+\mbox{Im}(\xi_1)s & \mbox{Re}(\xi_1\sigma_3)-\mbox{Re}(\xi_1)s \\
0 & & & \mbox{Im}(\xi_3\overline{\sigma}_3) -\mbox{Im}(\xi_3)s & \mbox{Re}(\xi_3\overline{\sigma}_3)-\mbox{Re}(\xi_3)s
\end{bmatrix}
\]
\end{example}
As for the original linearization, the pencil obtained by the removal of the first row, is not full rank.

A similar operation can be performed as in the standard case so that the size of the linearization is reduced by one and the last rows form a full rank dual linear basis.
The idea is to remove the last but first column, by a linear combination of the basis functions, so that $\psi_{d-1}$ is mapped to zero, using
\[
  \psi_{d-1} = 1 - \psi_1 - \psi_3 - \cdots - \psi_{d-3}.
\]
Using this linear combincation, we define a transformation so that
\[
  \begin{bmatrix} 1 & 0 \\ z & Z \end{bmatrix}  \begin{pmatrix} 1 \\ \psi_1 \\ \vdots \\ \psi_{d-2} \\ \psi_{d-1} \\ \psi_d \end{pmatrix} = \begin{pmatrix} 1 \\ \psi_1 \\ \vdots \\ \psi_{d-2} \\ 0 \\ \psi_d \end{pmatrix}.
\]
Multiplying ${L}_3$ on the right with the transformation performs a rank one update by subtracting the outer product of column $d$ and row $2$ from ${L}_3$.
Column $d$ has nonzero elements for rows one, two, and the last three rows.
Therefore, the outer product modifies the first, second and last three rows of ${L}_3$.
Denote $\mbox{Re}(r(\sigma_j))$ by $\rho_j$ and $\mbox{Im}(r(\sigma_{j})$ by $\eta_j$.
The first row is transformed from
\[
  \begin{matrix} 0 & \rho_1 & \eta_1 & \cdots & \rho_{d-3} & \eta_{d-3} & \rho_{d-1} & \eta_{d-1} \end{matrix}
\]
to
\[
  \begin{matrix} \rho_{d-1} & \rho_1-\rho_{d-1} & \eta_1 & \cdots & \rho_{d-3}-\rho_{d-1} & \eta_{d-3} &  \rho_{d-1} & \eta_{d-1}  \end{matrix},
\]
The second row is transformed from
\[
  \begin{matrix} -1 & 1 & 0 & \cdots & 1 & 0 & 1 & 0 \end{matrix}
\]
to
\[
  \begin{matrix} 0 & 0 & 0 & \cdots & 0 & 0 & 1 & 0 \end{matrix},
\]
i.e., only the value $1$ corresponding to $\psi_{d-1}$ remains.
%\[
%h_l(s) T^{-1}=\begin{bmatrix}h_{l,1}(s) & h_{l,2}(s) \end{bmatrix} = \begin{bmatrix} \mbox{Im}(\xi_l\overline{\sigma}_l) -\mbox{Im}(\xi_{l})s \mbox{Re}(\xi_{l}\overline{\sigma}_{l})-\mbox{Re}(\xi_{l})s \end{bmatrix}
%\]
%The last row is transformed from
%\[
%  \begin{matrix} 0 & 0 & 0 & \cdots & 0 & 0 & h_{d-1,1}(s) & h_{d-1,2}(s) \end{matrix}
%\]
%to
%\[
%  \begin{matrix} h_{d-1,1}(s) & -h_{d-1,1}(s) & 0 & \cdots & -h_{d-1,1}(s) & 0 & h_{d-1,1}(s) & h_{d-1,2}(s) \end{matrix}
%\]
Since $\psi_{d-1}$ is mapped to zero, we can remove the $d$th column.
This leaves us with a zero second row, which can also be removed.
\begin{lemma}
The dual linear basis consisting of the last $d-1$ rows of $L_4$ has full rank.
\end{lemma}
\begin{proof}
By moving the first column of ${L}_4$ in between columns of $d-1$ and $d$, we have the matrix
\[
  \begin{bmatrix} \star & \cdots & \cdots & \star & \star \\
    -\tilde{L}_{3,1} & \tilde{L}_{1,3} \\
    & \ddots & \ddots \\
          &        & \tilde{L}_{d-3,d-5} & \tilde{L}_{d-5,d-3} \\
    \star & \cdots & \star & \tilde{L}_{d-1,d-3} & \tilde{L}_{d-3,d-1} \\
    \star & \cdots & \star & & \tilde{h}_{d-1}^T
    \end{bmatrix},
\]
where $\star$ indicates nonzero blocks.
We have to inspect the submatrix obtained by removal of the first row.
By applying the inverse transformations that we used to form ${L}_3$ from ${L}_2$, we obtain the pencil
\[
  \begin{bmatrix}
    -L_{3,1} & L_{1,3} \\
    & \ddots & \ddots \\
    \star & \star & L_{d-1,d-3} & L_{d-3,d-1} \\
    \star & \star & & h_{d-1}^T
    \end{bmatrix}.
\]
The stars are a copy, possibly with a sign change, of the last but first column.
We can add a multiple of the last but first column to the columns to remove the stars.
This does not change the rank.
As a result, we have the matrix
\[
  \begin{bmatrix}
    -L_{3,1} & L_{1,3} \\
    & \ddots & \ddots \\
          &       & L_{d-1,d-3} & L_{d-3,d-1} \\
          &       & & h_{d-1}^T
    \end{bmatrix}.
\]
We are going to show that this matrix has full rank for all $s\in\mathbb{C}\cup\{\infty\}$.
The main diagonal blocks have full rank for
$s$ not in $\{\sigma_1,\ldots,\sigma_{d-1}\}$, and the off-diagonal blocks are full rank for $s\not\in\{\sigma_2,\ldots,\sigma_d\}$.
If $s=\sigma_j$, e.g., with $j\in\{1,\ldots,d-2\}$, then we use the fact that, after the elmination of column $j$, we still have a full rank matrix.
\end{proof}

\subsection{Real valued linearization for partial fraction basis}
Next to the barycentric rational basis, we also use the partial fraction basis for the AAA least squares method.
The AAA poles are the eigenvalues of the real pencil $C_1-s C_0$ and form a set of complex conjugate pairs.
We can use partial fractions such as
\[
  \Phi = \begin{pmatrix} 1 & \frac{\xi_1}{s-\sigma_1} & \cdots & \frac{\xi_d}{s-\sigma_d} \end{pmatrix}^T,
\]
where $\xi_1,\ldots,\xi_d$ are the positive nonzero real weights.
%An alternative choice could be the inverse of Newton polynomials, e.g.,
%\[
%\left\{ 1, \frac{1}{s - \sigma_1}, \frac{1}{s - \sigma_1}\frac{1}{s - \sigma_2}, \ldots,\prod_{i=1}^d\frac{1}{s-\sigma_i} \right\}.
%\]
%This choice would be useful when poles are multiple or consists of tight clusters.
%We do not explore this choice in the paper.
A rational function $r(s)$ is expressed as
\[
r(s) = \gamma_0 + \gamma_1\frac{\xi_1}{s-\sigma_1} + \cdots + \gamma_d\frac{\xi_d}{s-\sigma_d},
\]
where $\gamma_0,\ldots,\gamma_d$ are the coefficients of the function expressed in the basis.
A linear dual basis of $\Phi(s)$ is
\begin{equation}\label{eq:lin pf}
L_{\rm pf} = \begin{bmatrix}
  \xi_1 & \sigma_1 - s & 0 \\ \xi_2 & 0 & \sigma_2-s & 0 \\
  \vdots & & & \ddots \\
  \xi_d &  &  & & & \sigma_d -s
  \end{bmatrix}
\end{equation}
Since the poles form a complex conjugate set, we can form a real linear dual basis pencil, for a modified basis, as for the
barycentric formulation.
Let the poles be ordered so that $\sigma_1,\ldots,\sigma_r$ are real and $\sigma_{r+1},\ldots,\sigma_{d}$ are a set of complex
conjugate poles.
Define
\begin{eqnarray*}
  \mathbb{T}_R & = & \mbox{diag}\left( \underbrace{1,\ldots,1}_{\mbox{$r+1$ times}}, \underbrace{T,\ldots,T}_{\mbox{$d-r+1$ times}}\right) \\
  \mathbb{T}_L & = & \mbox{diag}\left( \underbrace{1,\ldots,1}_{\mbox{$r$ times}}, \underbrace{T,\ldots,T}_{\mbox{$d-r+1$ times}}\right) \\
\end{eqnarray*}
where $T$ is defined in \eqref{eq:T}.
We define the new basis of real functions
\[
{\Psi} = \mathbb{T}_R {\Phi}.
\]
As for the barycentric basis, the pencil $\mathbb{T}_L L_{\rm pf} \mathbb{T}_R^{-1}$ is real.
It is easy to see that \eqref{eq:lin pf} has full rank for all $s\in\mathbb{C}\cup\{\infty\}$, and so is the pencil after transformation.

For the inverse Newton polynomials
\[
  \Phi = \begin{pmatrix} 1 & \frac{\xi_1}{s-\sigma_1} & \frac{\xi_2}{(s-\sigma_1)(s-\sigma_2)} & \cdots & \frac{\xi_d}{(s-\sigma_1)\dots(s-\sigma_d)} \end{pmatrix}^T,
\]
we use the following linear dual basis
\begin{equation}\label{eq:lin in}
L_{\rm in} = \begin{bmatrix}
  \xi_1 & \sigma_1 - s & 0 \\ 0 & \xi_2 & \xi_1(\sigma_2-s) & 0 \\
  \vdots & & & \ddots \\
   &  &  & & \xi_d & \xi_{d-1}(\sigma_d -s)
  \end{bmatrix},
\end{equation}
so that $L_{\rm in}\Phi =0$.
If all poles are complex conjugate pairs, we have the basis
\[
  \Phi = \begin{pmatrix} 1 & \frac{\xi_1}{s-\sigma_1} & \frac{\xi_2}{(s-\sigma_1)(s-\sigma_2)} & \cdots & \frac{\xi_d}{(s-\sigma_1)\dots(s-\sigma_d)} &
        \frac{\xi_1}{s-\overline{\sigma}_1} & \frac{\xi_2}{(s-\overline{\sigma}_1)(s-\overline{\sigma}_2)} & \cdots & \frac{\xi_d}{(s-\overline{\sigma}_1)\dots(s-\overline{\sigma}_d)} \end{pmatrix}^T,
\]
and the following linear dual basis
\begin{equation*}
L_{\rm in} = \begin{bmatrix}
  \xi_1 & \sigma_1 - s & 0 & 0 \\
  \xi_1 &      0       & \overline{\sigma}_1 - s & 0 \\
    &   \xi_2  &       & \xi_1(\sigma_2-s) & 0 \\
    &      &   \xi_2   & 0 & \xi_1(\sigma_2-s) & 0 \\
      & & & \ddots  & & \ddots \\
  &  &  & & \xi_d & 0 & \xi_{d-1}(\sigma_d -s) \\
  &  &  & & 0 & \xi_d & 0 & \xi_{d-1}(\overline{\sigma}_d -s)
  \end{bmatrix}.
\end{equation*}
The basis is mapped to real functions and the linear dual basis to a real pencil as for the barycentric formulation.
For the example with weights one and poles $1$, $2\pm3\imath$, the real linear pencil is
\[
  \begin{bmatrix} 1 \\ & T \end{bmatrix} \begin{bmatrix}
     1 & 1 - s & 0 & 0 \\
     0 & 1     & 2+3\imath - s & 0 \\
     0 & 1     & 0  & 2-3\imath -s
  \end{bmatrix} \begin{bmatrix} 1 \\ & 1 \\ & & T^{-1} \end{bmatrix} =
  \begin{bmatrix}
     1 & 1 - s & 0 & 0 \\
     0 & 2     & 2 - s & 3 \\
     0 & 0     & -3  & 2 -s
  \end{bmatrix}.
\]
The basis of real rational functions is the set of real functions
\[
  1,\quad 1/(s-1), \quad 2(s-2)/(s-1)((s-2)^2+9), \quad -6/(s-1)((s-2)^2+9).
\]
As for the partial fraction formulation, it is not hard to see that the linear dual basis is a pencil of full rank.

\subsection{Extended AAA and Filtered AAA}
Reconsider the rational function \eqref{eq:schur complement}.
Assume that $C_1$ is nonsingular, otherwise shift with XXXXX
If the dual linear basis has full rank, then $C_1$ is invertible (NOT CORRECT IF $c_1\neq0$), so, we have
\begin{eqnarray*}
\begin{pmatrix} s r(s) \\ 0 \end{pmatrix}
& = & \begin{bmatrix} \alpha s & a^T s \\ c_0+s c_1 & C_0 + s C_1 \end{bmatrix} \begin{pmatrix} 1 \\ \Phi(s) \end{pmatrix} \\
& = & \begin{bmatrix} \alpha s - a^T C_1^{-1} (c_0+s c_1) & -a^T C_1^{-1} C_0 \\ c_0+sc_1 & C_0 + s C_1 \end{bmatrix} \begin{pmatrix} 1 \\ \Phi(s) \end{pmatrix},
\end{eqnarray*}
which is obtained by adding a linear combination of the two block rows of the equations to the first block row.
In words, $s r(s)$ can be written as a linear combination of the basis functions $\Phi(s)$ plus a linear term in $s$.
Also,
\begin{eqnarray*}
\begin{pmatrix} s^2 r(s) \\ 0 \end{pmatrix}
& = & \begin{bmatrix} \alpha s^2 - a^T C_1^{-1} (c_0+sc_1) s & -a^T C_1^{-1} C_0 s \\ c_0+sc_1 & C_0 + s C_1 \end{bmatrix} \begin{pmatrix} 1 \\ \Phi(s) \end{pmatrix}, \\
& = & \begin{bmatrix} \alpha s^2 - a^T C_1^{-1} (c_0+sc_1) s + a^T C_1^{-1} C_0 C_1^{-1} (c_0+sc_1)  & a^T C_1^{-1} C_0 C_1^{-1} C_0 \\ c_0+sc_1 & C_0 + s C_1 \end{bmatrix} \begin{pmatrix} 1 \\ \Phi(s) \end{pmatrix}.
\end{eqnarray*}
In words, $r(s)$, $sr(s)$ and $s^2 r(s)$ are linear combinations of $1$, $s$, $s^2$ and $\phi_1,\ldots,\phi_d$.
That is, for $r$, $sr$ and $s^2r$ there are $\tilde{a}$, $\alpha_0$, $\alpha_1$ and $\alpha_2$, respectively, so that
\[
\begin{bmatrix}
\alpha_0 + s \alpha_1 + s^2 \alpha_2 & a^T \\ c_0+sc_1 & C_0 + s C_1
\end{bmatrix}
\]
represents $r$, $sr$ and $s^2r$, respectively.

These observations give more possibilities for approximating nonlinear matrices, when taking into account a degree two matrix polynomial part in the approximation.
Quadratic matrices are the usual representation of a vibration problem using FEM, so, it makes sense to add a polynomial of degree two in the approximation.
In numerical experiments, we noticed that AAA does not always make a good approximation, with stable poles, for the nonlinear terms.
In Section~\ref{sec:numex}, we show a case with a nonlinear function with a factor $s^2$ that is not approximated well by a rational function with stable poles.
Therefore, we extend the barycentric form by adding a polynomial term that captures the linear and quadratic behaviour of the functions.

To achieve this, once, $c_0$, $c_1$, $C_0$, and $C_1$ are computed from Set Valued AAA or Weighted AAA, possibly with a filtering of the unstable eigenvalues,
we solve the least squares problem
\[
\min \sum_{s\in\mathbf{Z}} \left|\alpha_0 + s \alpha_1 + s^2 \alpha_2 + a^T \Phi(s) - g(s)\right|^2.
\]
This problem can be solved for the barycentric form or, in the case of unstable poles, the partial fraction form or inverted Newton polynomial form, after elimination of the unstable poles.
After the approximation, we obtain the following rational matrix:
\begin{eqnarray*}
\tilde{A}_0 + s \tilde{A}_1 + s^2 \tilde{A}_2 + \sum_{j=1}^m A_{-j} a_j^T \Phi(s)
\end{eqnarray*}
with
\begin{eqnarray*}
\tilde{A}_k & = & A_k + \sum_{j=1}^m \alpha_k^{(j)} A_{-j} \quad,\quad k=0,1,2,
\end{eqnarray*}
and
$\alpha_0^{(j)}$, $\alpha_1^{(j)}$, $\alpha_2^{(j)}$ and $a_j$ the coefficients for the approximation of $g_j$.
The solution of the least squared problem does not modify $C_0$ and $C_1$.
Therefore, the poles of the rational approximation are not modified.
As a consequence, stability of the poles is preserved.

\section{Time integration}\label{sec:time}

We show in this section that the solution of the linear system of ODEs is connected to \eqref{eq:sys}.
Theorem~\ref{th:connection} shows the connection between the nonlinear frequency dependency in \eqref{eq:sys} and the linear initial value problem
\begin{equation}\label{eq:lin ode}
-\mathbf{E} \frac{d \mathbf{x}}{dt} + \mathbf{A} \mathbf{x} = \mathbf{b}(t)
  \quad,\quad \mathbf{x}=\mathbf{x}_0,
\end{equation}
for a particular choice of initial values, with $\mathbf{A}$ and $\mathbf{E}$ defined in \eqref{eq:rat lin}.
The proof of the theorem relies on the following two lemmas, that consider two specific situations.
The first is the situation of a periodic solution at a given angular frequency $\omega=2\pi f$, the second is the case of zero initial values and its derivative.

The size of the linearization is a multiple of the size of the original problem.
This has an impact on the cost of the simulation of the linear system.
Fortunately, there is a large amount of structure in the linearization pencil, which can be exploited in numerical computations.
In a time integration method such as backward Euler or Crank-Nicholson, the main operation has the form
\[
  \mathbf{w} = (h \mathbf{A} + \mathbf{E})^{-1}(\mathbf{b} + \mathbf{E} \mathbf{v}).
\]
The efficient implementation relies on the UL factorization of $h \mathbf{A} + \mathbf{E}$ \cite{bemm15}.
Let
\[
  \check{A} = a_1^T \otimes A_{-1} + \cdots + a_m^T \otimes A_{-m}.
\]
Then, the linearization is factorized as
\begin{multline}\label{eq:UL factorization}
\left[\begin{array}{c|cc} \tilde{A}_0+s A_1 & s A_2 & \check{A} \\\hline
                           s I & -I & 0 \\
                           (c_0+sc_1)\otimes I &0 & (C_0+sC_1)\otimes I
                \end{array}\right] \\  = \underbrace{\left[\begin{array}{c|cc}
                         I & -s A_2 & \check{A} ((C_0+sC_1)^{-1}\otimes I) \\\hline
                         0 & I
                \end{array}\right]}_{\mathbf{U}}
                \left[\begin{array}{c|cc}
                         R(s) & 0 & 0 \\\hline
                           s I & -I & 0  \\
                         (c_0+sc_1)\otimes I &0 & (C_0+sC_1)\otimes I
                \end{array}\right].
\end{multline}
The computation of $\mathbf{w}$ can use the UL factorization.
We deduce that it requires linear combinations of blocks of $\mathbf{v}$ and $\mathbf{b}$,
matrix vector products with the coefficient matrices, and the solution of a linear system with $R(s)$ with $s=-h^{-1}$.
That is, the cost is much lower than the cost that one would expect from a linear solver directly applied to $h \mathbf{A} + \mathbf{E}$.

%\begin{lemma}
%Let
%\[
%  \mathbf{x} = \begin{pmatrix} x_0 \\ s x_0 \\ \mathbf{x}_1(s) \end{pmatrix}
%\]
%such that
%\[
%  (\mathbf{A} - s \mathbf{E}) \mathbf{x} = e_1 \otimes A(s) x_0.
%\]
%Then
%\[
%  \mathbf{E} \mathbf{x} = 
%\]
%\end{lemma}
%\begin{proof}
%We have that the first block row of $\mathbf{E} \mathbf{x}$ is
%\[
%  -\tilde{A}_1 x_0 - s A_2 x_0 + E_{-1} r_1(s) x_0 + \cdots E_{-m} r_m(s) x_0.
%\]
%\end{proof}

Define
\[
  \mathbf{\Phi} = \begin{pmatrix} 1 \\ s \\ \Phi(s) \end{pmatrix}:\mathbb{C}\to\mathbb{C}^{\tilde{d}+2}
\]
satisfying $(C_0+sC_1)\Phi(s) + (c_0+sc_1)=0$ for all $s\in\mathbb{C}$.
\begin{lemma}\label{le:case1}
Consider the initial value problem
\begin{eqnarray*}
  \mathbf{A} \mathbf{x} -\mathbf{E} \frac{d \mathbf{x}}{dt} & = & e_1\otimes b e^{\imath\omega t} \\
    \mathbf{x}(0) & = & \mathbf{\Phi}(\imath\omega) \otimes \hat{x},
\end{eqnarray*}
where $\hat{x}$ is the solution of
\[
R(\imath\omega) \hat{x} = b.
\]
Then, the solution of the initial value problem
is $\mathbf\Phi(\imath\omega) \otimes \hat{x}e^{\imath\omega t}$.
The Laplace transform $\mathcal{L}(x)$ satisfies
\begin{equation}\label{eq:lapl-harmonic}
R(s) \mathcal{L}(x) = \mathcal{L}(be^{\imath\omega t}) + A_1 x(0) + A_2(sx(0)+\dot{x}(0)) - \check{A}( ((C_0+sC_1)^{-1}\otimes I) \Phi(\imath\omega)\otimes x(0)).
\end{equation}
\end{lemma}
\begin{proof}
The initial value problem
\[
  \mathbf{A} \mathbf{x} -\mathbf{E} \frac{d \mathbf{x}}{dt} = e_1\otimes b e^{\imath\omega t}
\]
has solution $\mathbf{x}(t) = (\mathbf{A}-\imath\omega\mathbf{E})^{-1} (e_1\otimes b) e^{\imath\omega t}$
if the initial values are $(\mathbf{A}-\imath\omega\mathbf{E})^{-1} (e_1\otimes b)$.
Following the structure of $\mathbf{A}$ and $\mathbf{E}$, we have that
$\mathbf{x}(t) = (\mathbf{\Phi}(\imath\omega) \otimes \hat{x}) e^{\imath\omega t}$.

The Laplace transform of the initial value problem is
\[
  (\mathbf{A} - s \mathbf{E}) (\mathbf{\Phi}(\imath\omega)\otimes \hat{x})/(s-\imath\omega) = e_1\otimes b / (s-\imath\omega) - \mathbf{E} (\mathbf{\Phi}(\imath\omega)\otimes \hat{x}).
\]
After the multiplication with $\mathbf{U}^{-1}$ on the left,
The first block row of this equation leads to \eqref{eq:lapl-harmonic}.
This proves the lemma.
\end{proof}

\begin{lemma}\label{le:case2}
The initial value problem
\begin{eqnarray}\label{eq:system-theorem1}
  \mathbf{A} \mathbf{x} -\mathbf{E} \frac{d \mathbf{x}}{dt} & = & \mathbf{b}(t) = e_1\otimes b(t) \\
    \mathbf{x}(0) & = & 0
\end{eqnarray}
has as solution
\[
  \mathbf{x}(t) = \begin{pmatrix} x(t) \\ \dot{x}(t) \\ \mathbf{x}_1(t) \end{pmatrix}
\]
where $\mathcal{L}(x(t))$ satisfies the algebraic equation
\begin{equation}\label{eq:system-laplace-theorem}
  R(s) \mathcal{L}(x) = \mathcal{L}(b).
\end{equation}
\end{lemma}
\begin{proof}
Again decompose
\[
  \mathbf{x} = \begin{pmatrix} x \\ \frac{dx}{dt} \\ \mathbf{x}_1 \end{pmatrix}.
\]
Taking into account the structure of $\mathbf{A}$ and $\mathbf{E}$, we have that
\begin{eqnarray}
  (C_0\otimes I) \mathbf{x}_1 + (C_1\otimes I) \frac{d\mathbf{x}_1}{dt} & = & -c_0 \otimes x - c_1 \otimes \frac{dx}{dt}, \label{eq:x_1 ode}\\
    \mathbf{x}_1(0) & = & 0. \nonumber
\end{eqnarray}
The Laplace transform of $\mathbf{x}_1$ is therefore $\mathcal{L}(\mathbf{x}_1) = -(C_0+sC_1)^{-1} (c_0+sc_1) \otimes \mathcal{L}(x)$.

With the given structure of $\mathbf{x}(t)$, the Laplace transform of the
first block row of \eqref{eq:system-theorem1} produces \eqref{eq:system-laplace-theorem}.
\end{proof}

\begin{theorem}\label{th:connection}
Consider the initial value problem
\begin{eqnarray}\label{eq:system-theorem}
  \mathbf{A} \mathbf{x} -\mathbf{E} \frac{d \mathbf{x}}{dt} & = & \mathbf{b}(t) = e_1\otimes b(t) \\
    \mathbf{x}(0) & = & \begin{pmatrix} x(0) \\ \dot{x}(0) \\ \mathbf{x}_1(0) \end{pmatrix} \nonumber\\
  \mathbf{x}_1(0) & = & \mathcal{F}^{-1}( \Phi(\imath\omega)\otimes \mathcal{F}( \widetilde{x}(t) ) )_{t=0} \nonumber,
  %\mathbf{x}_1(t) & = & \frac{1}{\pi} \int_{-\infty}^\infty  \Phi(\imath\omega)\otimes \mbox{Re}(\int_{-\infty}^0 x(t) e^{-\imath\omega t} dt \nonumber\\
\end{eqnarray}
where $\mathcal{F}$ is the Fourier transform, and $\widetilde{x}$ is a function that is equal to $x(t)$ for $t\leq 0$.
Let the solution $\mathbf{x}(t)$ have the form
\[
  \mathbf{x}(t) = \begin{pmatrix} x(t) \\ \dot{x}(t) \\ \mathbf{x}_1(t) \end{pmatrix}.
\]

If the Fourier and inverse Fourier transforms that define $\mathbf{x}_1(0)$ exist and converge,
we have the following algebraic equation for the Laplace transform of $x(t)$:
\[
R(s) \mathcal{L}(x) = \mathcal{L}(b(t)) + \tilde{A}_1 x(0) + \tilde{A}_2(sx(0)+\dot{x}(0)) - \check{A}( ((C_0+sC_1)^{-1}\otimes I) \mathbf{x}_1(0) ).
\]
\end{theorem}
Note that $\mathbf{x}_1$ is the solution of \eqref{eq:x_1 ode} integrated from $-\infty$ to $0$.
It is a linear combination of purely harmonic terms as defined in Lemma~\ref{le:case1}.
In practice, one could start from $-T$ with $x(t)=0$ for $t\leq-T$ and $T$ large enough so that $\|x(-T)\|$ is small.

In the expression of the Laplace transform, the terms in the right-hand correspond to what we would obtain from a second order equation
plus a term from the rational approximation.
The terms from the rational approximation can also be written as
\[
    \check{A}( ((C_0+sC_1)^{-1}\otimes I) \mathbf{x}_1(0) ) = \sum_{i=1}^m ((a_i^T (C_0+sC_1)^{-1}) \otimes A_i ) \mathbf{x}_1(0).
\]

\noindent
\begin{proof}
%If $x(t) = \hat{x}(\imath\omega) e^{\imath\omega t}$ for $t\leq 0$, then the lower part of the initial values vector is
%\[
%  \mathbf{x}_1(0) = \Phi(\imath\omega)\otimes \hat{x}(\imath\omega).
%\]
%
Let $\widetilde{x}(t)$ be an extension of $x(t)$ from $t\in(-\infty,0]$ to $\mathbb{R}$.
The Fourier transform of $\widetilde{x}(t)$ can be used to represent $x(t)$ as a sum (or integral) of harmonic functions for $t<0$.
Let $\check{x}(\omega)=\mathcal{F}(\widetilde{x})$ be the Fourier transform of $\widetilde{x}$, then $\widetilde{x}(t)=\mathcal{F}^{-1}(\check{x}(\omega))$.

Following Lemma~\ref{le:case1}, we have that
\begin{eqnarray*}
\mathbf{A} \mathbf{h} -\mathbf{E} \frac{d \mathbf{h}}{dt} & = & e_1\otimes R(\imath\omega)\check{x}(\omega) e^{\imath\omega t} \label{eq:proof-theorem-1},\\
\mathbf{h}(0) & = & \mathbf{\Phi}(\imath\omega)\otimes \check{x}(\omega), \nonumber
\end{eqnarray*}
with solution $\mathbf{h}(t) = \mathbf{\Phi}(\imath\omega) \otimes \check{x}(\omega) e^{\imath\omega t}$.
We now make a linear combination of this equation for $\omega$ varying from $-\infty$ to $\infty$.
Consider
\begin{eqnarray*}
\mathbf{A} \widetilde{\mathbf{h}} -\mathbf{E} \frac{d \widetilde{\mathbf{h}}}{dt} & = & e_1\otimes \mathcal{F}^{-1}( R(\imath\omega) \check{x}(\omega) ) \\
\widetilde{\mathbf{h}}(0) & = & \mathcal{F}^{-1}(\mathbf{\Phi}(\imath\omega)\otimes \check{x}(\omega))|_{t=0} \nonumber
\end{eqnarray*}
with solution $\widetilde{\mathbf{h}}(t) = \mathcal{F}^{-1}(\mathbf{\Phi}(\imath\omega)\otimes \check{x}(\omega))$.

Now subtract $\tilde{b}=\mathcal{F}^{-1}( R(\imath\omega) \check{x}(\omega) )$ from the right hand side of \eqref{eq:system-theorem}:
\begin{eqnarray*}
\mathbf{A} \mathbf{z} -\mathbf{E} \frac{d \mathbf{z}}{dt} & = & e_1\otimes (b(t)-\tilde{b}(t)) \\
\mathbf{z}(0) & = & 0.
\end{eqnarray*}
The solution of \eqref{eq:system-theorem} is then $x = \widetilde{h} + z$.

For the Laplace transform of the equation, we take into account that the Laplace transform is a linear operator and therefore:
\begin{eqnarray*}
\mathbf{A}( \mathcal{L}(\widetilde{h}) + \mathcal{L}(z) ) - s \mathbf{E} ( \mathcal{L}(\widetilde{h}) + \mathcal{L}(z) ) ) = \mathcal{L}(b-\tilde{b}(t))  + \mathcal{L}(\tilde{b}(t)) - \mathbf{E} \mathbf{\widetilde{h}}(0)
\end{eqnarray*}
The proves follows from Lemma~\ref{le:case1} applied to the superposition of the harmonic parts.
\end{proof}

The Laplace transform leads to a simple connection between the nonlinear frequency dependent model and a linear system of ODEs.
Note that the initial values associated with the rational terms use information from the past, which may not be available.
At time $t=0$, we can impose $x$ and its derivative, but we usually assume that the solution is zero for $t<0$, so that $\mathbf{x}_1(0)=0$.
If this is not the case, a practical way to compute $\mathbf{x}_1(0)$ is to use the differential equation \eqref{eq:x_1 ode},
with initial value zero at some time $T\ll 0$ with $T$ well enough in the past so that $|x(t)|$ for $t\leq T$ is small or zero.

Let us, as an example, assume that the system is at rest before $t=0$, with $x(t)=x_0$ and $\dot{x}(t)=0$ for $t\leq0$.
Then, we can use Lemma~\ref{le:case1} with $\omega=0$, i.e.,
\[
  \mathbf{x}_1(0) = \Phi(0) \otimes x_0.
\]

%The state vector of the linearization has a special structure that we now discuss.
%We assume that the right-hand side always has the form $e_1\otimes b(t)$ and that the initial values have the form \eqref{eq:system-theorem}.
%Then the solution is
%\[
%\]
%\begin{theorem}
%Let $Q\in\mathbb{R}^{n\times r}$ be such that $x(t) = Q \tilde{x}(t)$ for all $t$.
%Then,
%\[
%  \mathbf{x}_1(t) = (I_d\otimes Q) z(t)
%\]
%with $z(t)$ the solution of
%\begin{eqnarray*}
%  (C\otimes I_r) z - (D\otimes I_r) \frac{d z}{dt} & = & e_1 \otimes \tilde{x}(t) \\
%    z(0) & = & \int_{-\infty}^0 (\exp(-s D^{-1}C) D^{-1} e_1) \otimes \tilde{x}(s) ds  \end{pmatrix}\nonumber.
%\end{eqnarray*}
%\end{theorem}
%\begin{proof}
%The proof immediately follows from Theorem~\ref{th:connection}, where $x(t)$ is replaced by $Q\tilde{x}(t)$.
%\end{proof}

\section{Numerical examples}\label{sec:numex}

The numerical examples are run on a Macbook Pro 2,8 with GHz Quad-Core Intel Core i7 processor and 16 GB RAM of 2133 MHz, using the C++ library CORK++ \cite{CORKpp}, which includes the numerical examples source codes
of this paper.

We compare the following methods:
\begin{itemize}
\item AAA: this is weighted AAA for a set of functions;
\item AAA-LS: this is weighted AAA, followed by a least squares step on the barycentric formulation;
\item E-AAA: this is extended AAA-LS, i.e., weighted AAA followed by a least squares step for the barycentric formulation combined with a polynomial of degree two;
\item F-AAA: this is AAA followed by a least squares approximation on the partial fraction formulation using only the stable poles;
\item S-AAA: this is AAA followed by a least squares approximation on the partial fraction formulation where all unstable poles are flipped along the imaginary axis to the left half plane;
\item E-F-AAA: this is F-AAA, but now with a least squares approximation combined with a polynomial of degree two.
\item E-S-AAA: this is S-AAA, but now with a least squares approximation combined with a polynomial of degree two.
\end{itemize}
For F-AAA and E-F-AAA, we noticed no difference in accuracy between partial fractions and inverted Newton polynomials.
The numerical results shown are for the partial fraction form.
The difference between S-AAA and F-AAA lies in that F-AAA throws away unstable poles and reduces the dimension of the linearization, where S-AAA changes the sign of the real parts of the unstable poles.

\subsection{Sandwich beam}
\begin{figure}
\begin{center}
\includegraphics[width=0.7\textwidth]{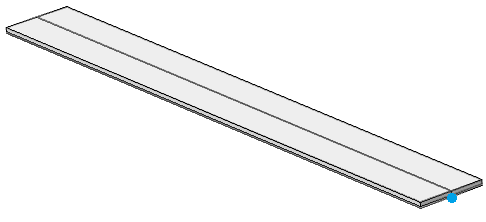}
\end{center}
\caption{Configuration of the beam structure}\label{fig:beam}
\end{figure}

The model consists of a clamped, thin, flat aluminum beam structure, consisting of two steel layers surrounding a damping layer, represented by a fractional derivative model \cite{Bagley1983}; see Figure~\ref{fig:beam}.
The model is discretized by the finite element method and is described by the algebraic equation
\begin{equation*}
A(\omega) = A_0 - \omega^2 A_2 + \frac{G_0 + G_{\infty} (i \omega \tau)^{\alpha}}{1 + (i \omega \tau)^{\alpha}} A_{-1},
\end{equation*}
with $A_0,A_1,A_2 \in \mathbb{R}^{168 \times 168}$ symmetric positive semi-definite matrices; see~\cite{vanbeeumen_meerbergen}.
Parameters are the static shear modulus $G_0 = 350.4$kPa, the asymptotic shear modulus $G_\infty = 3.062$MPa, the relaxation time $\tau = 8.23$ns and a fractional parameter $\alpha=0.675$.
The model is valid for frequency range $[200,30000]$Hz.

\begin{figure}
\begin{center}
\begin{tabular}{cc}
\begin{tikzpicture}[scale=0.8]
\begin{semilogxaxis}[xmin=0,xmax=5000,xlabel={$f$ (Hz)},ylabel=$|g_1(\imath 2\pi f)|$]
\addplot[mark=none,blue,thick] table[y expr=sqrt(\thisrow{real}*\thisrow{real}+\thisrow{imag}*\thisrow{imag})] {beam.dat} ;
\end{semilogxaxis}
\end{tikzpicture} &
\begin{tikzpicture}[scale=0.8]
\begin{semilogxaxis}[xmin=0,xmax=10000,xlabel={$f$ (Hz)},ylabel=$\phi(g_1(\imath 2\pi f))$]
\addplot[mark=none,blue,thick] table[y expr=rad(atan(\thisrow{imag}/\thisrow{real}))] {beam.dat} ;
\end{semilogxaxis}
\end{tikzpicture} \\
  modulus & phase
\end{tabular}
\end{center}
\caption{Modulus $|g_1|$ and phase $\phi(g_1)$ of the sandwich beam for $s=\imath2\pi f$, $f\in[10,10^4]$Hz.}\label{fig:f_sandwich_beam}
\end{figure}
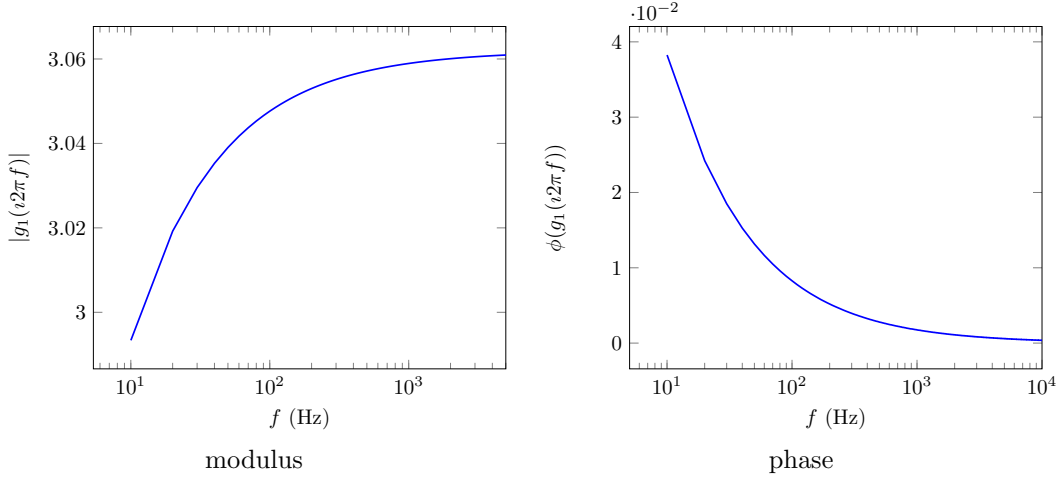

We show the rational approximation for the different choices of sample points.
We have to set up AAA so that the nonlinear functions are well approximated on the imaginary axis.
In principal, one could select a very large number of points on the imaginary axis, but we noticed that this leads to high computational costs.
The function $g_j$ typically varies slowly for large values of $\omega$; see, e.g., the
function
\[
  g_1(s) = \frac{G_0 + G_\infty(s\tau)^\alpha}{1+(s\tau)^\alpha}
\]
shown in Figure~\ref{fig:f_sandwich_beam}.
As a result, there is no need in using a high density of points for the higher frequencies.
We therefore use a logarithmic scale in the selection of points.
For the discretization in $n$ points of the interval $[\imath \omega_{\min},\imath \omega_{\max}]$ on the imaginary axis,
we use the set $Z_i = \imath 10^{\tilde{Z}_i}$, where $\{\tilde{Z}_i\}_{i=1}^M$ are an equidistant set on the
interval $[\log_{10}(\omega_{\min}),\log_{10}(\omega_{\max})]$, with $\omega=2\pi f$.
%Function $g_1$ has a singularity at $s=0$ and a pole at
%\[
%  s=\cos(\pi/\alpha)/\tau+\imath\sin(pi/\alpha)/\tau\simeq -0.0070650 - 0.12130\imath,
%\]
%which has a negative real part.
Note that $s=0$ is excluded from the domain of AAA, so that the approximation is not corrupted by the singularity at the origin.
%\begin{figure}
%\begin{center}
%\begin{tikzpicture}
%\begin{axis}[axis x line=middle, axis y line=right, xmin=-11,xmax=3, xlabel={$\mbox{Re}(\sigma_j)$}, ylabel={$\mbox{Im}(\sigma_j)$},
%                xtick={-10,-8,-6,-4,-2,0,2}, xticklabels={{\footnotesize$-10^{10}$},{\footnotesize$-10^{8}$},{\footnotesize$-10^{6}$},{\footnotesize$-10^{4}$},{\footnotesize$-10^{2}$},{\footnotesize$-10^{0}$},{\footnotesize$-10^{-2}$}}]
%\addplot[color=red,mark=*,only marks] table[x=xl,y=y] {beam_poles.dat} ;
%\end{axis}
%\end{tikzpicture}
%\end{center}
%\caption{Poles of the rational approximation for the damping function of the sandwich beam}\label{fig:beam_poles}
%\end{figure}

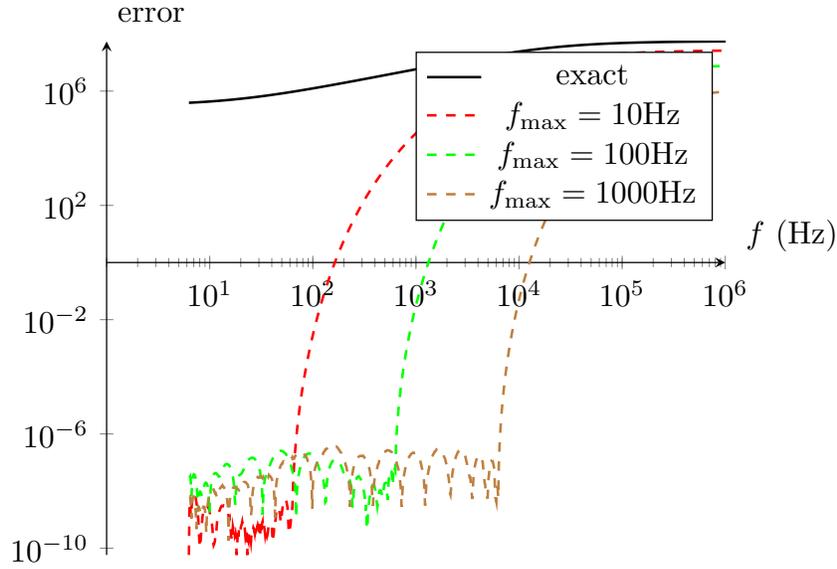
\begin{figure}
\begin{center}
\begin{tabular}{cc}
\begin{tikzpicture}[scale=1.2]
\begin{loglogaxis}[axis lines=center, xlabel={$f$ (Hz)}, ylabel={error},xmin=1,xmax=1000000]
\addlegendentry{exact}
\addplot[color=black,thick,ticks=none] table[x=s,y=n,col sep=comma] {beam_aaa_diff_10.dat} ;
\addlegendentry{$f_{\max}=10$Hz}
\addplot[color=red,thick,dashed,ticks=none] table[x=s,y=e,col sep=comma] {beam_aaa_diff_10.dat} ;
\addlegendentry{$f_{\max}=100$Hz}
\addplot[color=green,thick,dashed,ticks=none] table[x=s,y=e,col sep=comma] {beam_aaa_diff_100.dat} ;
\addlegendentry{$f_{\max}=1000$Hz}
\addplot[color=brown,thick,dashed,ticks=none] table[x=s,y=e,col sep=comma] {beam_aaa_diff_1000.dat} ;
\end{loglogaxis}
\end{tikzpicture}
\end{tabular}
\end{center}
\caption{Error of AAA approximation for different values of $f_{\max}$}\label{fig:beam-aaa}
\end{figure}

\paragraph{AAA approximation}
We used $n_Z=100$, $1000$, or $10000$ sample points in a logarithmic scale for intervals $[\imath f_{\min},\imath f_{\max}]$ with $f_{\min}=1$Hz and $f_{\max}=10,10^2,10^3,10^4$Hz.
The approximation for higher frequencies is good, even when only sample points are chosen in the lower spectrum.
This is illustrated by Figure~\ref{fig:beam-aaa}, that shows the error for two intervals on the frequency axis.
Note that with a tolerance $\tau=0.0$, the AAA algorithm stops when Froissart doublets are found. In this case, we had an unstable linearization for $f_{\max}=10$ and $n_Z=100$.

We used initial values zero for $x$ and its derivative and right-hand side $b(t)=b_0\cos(2\pi f t)$ with $f=10$Hz and $b_0$ zero everywhere except for element $470$, which has value $0.7071067811865475$.
We used the Crank-Nicholson method with time step $10^{-3}$.
Table~\ref{tab:stability} compares different choices for the AAA approximation, where we used the relative tolerance $\tau=10^{-13}$ as stopping criterion.
\begin{table}
\caption{Size of AAA approximation and stability for various choices of maximum frequency $f_{\max}$ and number of
  sample points $n_Z$}\label{tab:stability}
\begin{center}
\begin{tabular}{crcc}
$f_{\max}$ & \multicolumn{1}{c}{$n_Z$} & $d$ & Stable?\\\hline
$10^1$ & 100 & 18 & Yes \\
$10^1$ & 1000 & 18 & Yes \\
$10^1$ & 10000 & 18 & Yes \\
$10^2$ & 100 & 24 & Yes \\
$10^2$ & 1000 & 24 & Yes \\
$10^2$ & 10000 & 24 & Yes \\
\end{tabular}\qquad
\begin{tabular}{crcc}
$f_{\max}$ & \multicolumn{1}{c}{$n_Z$} & $d$ & Stable?\\\hline
$10^3$ & 1000 & 32 & Yes \\
$10^3$ & 10000 & 32 & Yes \\
$10^4$ & 1000 & 38 & Yes \\
$10^4$ & 10000 & 41 & Yes \\
\end{tabular}
\end{center}
\end{table}

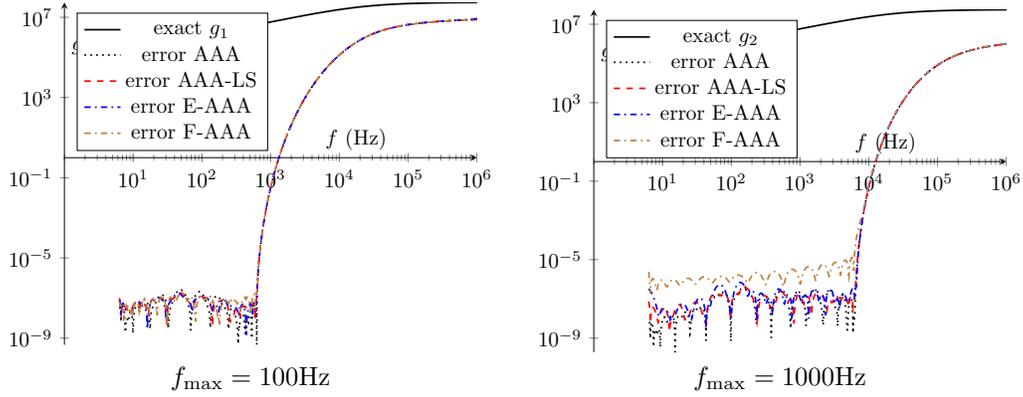
\begin{figure}
\begin{center}
\begin{tabular}{cc}
\begin{tikzpicture}[scale=0.8]
\begin{loglogaxis}[axis lines=center, xlabel={$f$ (Hz)}, ylabel={$g_1$},xmin=1,xmax=1000000,legend pos=north west]
\addlegendentry{exact $g_1$}
\addplot[color=black,thick,ticks=none] table[x=s,y=n,col sep=comma] {beam2/aaa_diff_100_10000.dat} ;
\addlegendentry{error AAA}
\addplot[color=black,dotted,thick,ticks=none] table[x=s,y=e,col sep=comma] {beam2/aaa_diff_100_10000.dat} ;
\addlegendentry{error AAA-LS}
\addplot[color=red,thick,dashed,ticks=none] table[x=s,y=e,col sep=comma] {beam2/aaa_ls_diff_100_10000.dat} ;
\addlegendentry{error E-AAA}
\addplot[color=blue,thick,dashdotted,ticks=none] table[x=s,y=e,col sep=comma] {beam2/e_aaa_diff_100_10000.dat} ;
\addlegendentry{error F-AAA}
\addplot[color=brown,thick,dashdotted,ticks=none] table[x=s,y=e,col sep=comma] {beam2/f_aaa_diff_100_10000.dat} ;
\end{loglogaxis}
\end{tikzpicture} &
\begin{tikzpicture}[scale=0.8]
\begin{loglogaxis}[axis lines=center, xlabel={$f$ (Hz)}, ylabel={$g_2$},xmin=1,xmax=1000000,legend pos=north west]
\addlegendentry{exact $g_2$}
\addplot[color=black,thick,ticks=none] table[x=s,y=n,col sep=comma] {beam2/aaa_diff_1000_10000.dat} ;
\addlegendentry{error AAA}
\addplot[color=black,thick,dotted,ticks=none] table[x=s,y=e,col sep=comma] {beam2/aaa_diff_1000_10000.dat} ;
\addlegendentry{error AAA-LS}
\addplot[color=red,thick,dashed,ticks=none] table[x=s,y=e,col sep=comma] {beam2/aaa_ls_diff_1000_10000.dat} ;
\addlegendentry{error E-AAA}
\addplot[color=blue,thick,dashdotted,ticks=none] table[x=s,y=e,col sep=comma] {beam2/e_aaa_diff_1000_10000.dat} ;
\addlegendentry{error F-AAA}
\addplot[color=brown,thick,dashdotted,ticks=none] table[x=s,y=e,col sep=comma] {beam2/f_aaa_diff_1000_10000.dat} ;
\end{loglogaxis}
\end{tikzpicture} \\
    $f_{\max}=100$Hz &
    $f_{\max}=1000$Hz
\end{tabular}
\end{center}
\caption{Comparison of AAA, AAA-LS and E-AAA for $f_{\max}=100$ and $f_{\max}=1000$ and $N_Z=1000$.}\label{fig:beam-aaa-ls}
\end{figure}

As can be seen from Figure~\ref{fig:beam-aaa-ls}, there is no difference in quality between the methods.
Therefore, further numerical experiments use AAA.
F-AAA does not remove poles since they are all stable, but we have added this comparison because
F-AAA uses partial fractions, where the other methods use the barycentric form of the rational function.
We notice that the error is a factor 1000 larger for the lower frequencies than for the barycentric form.
This was also observed for other examples.

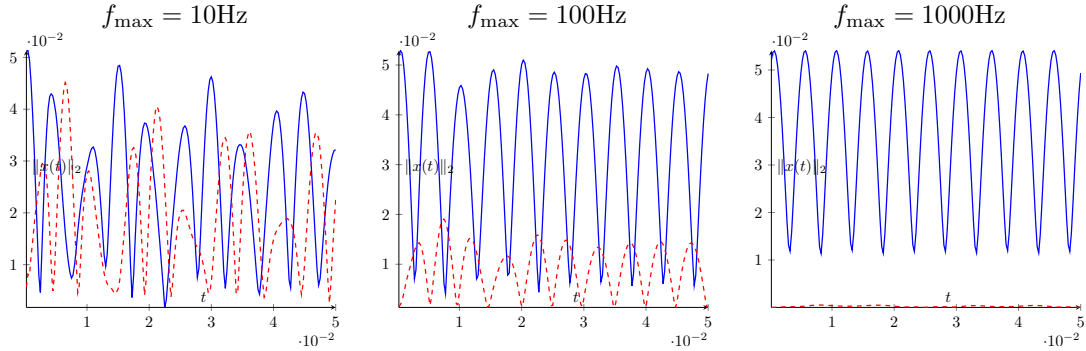
\begin{figure}
\begin{center}
\begin{tabular}{ccc}
$f_{\max}=10$Hz & $f_{\max}=100$Hz & $f_{\max}=1000$Hz \\
\begin{tikzpicture}[scale=0.6]
\begin{axis}[axis lines=center, xlabel={$t$}, ylabel={$\|x(t)\|_2$}]
\addplot[color=blue,thick,ticks=none] table[x=t,y=n,col sep=comma] {beam_periodic_10_1000_100.dat} ;
\addplot[color=red,thick,dashed,ticks=none] table[x=t,y=e,col sep=comma] {beam_periodic_10_1000_100.dat} ;
\end{axis}
\end{tikzpicture}
&
\begin{tikzpicture}[scale=0.6]
\begin{axis}[axis lines=center, xlabel={$t$}, ylabel={$\|x(t)\|_2$}]
\addplot[color=blue,thick,ticks=none] table[x=t,y=n,col sep=comma] {beam_periodic_100_1000_100.dat} ;
\addplot[color=red,thick,dashed,ticks=none] table[x=t,y=e,col sep=comma] {beam_periodic_100_1000_100.dat} ;
\end{axis}
\end{tikzpicture}
&
\begin{tikzpicture}[scale=0.6]
\begin{axis}[axis lines=center, xlabel={$t$}, ylabel={$\|x(t)\|_2$}]
\addplot[color=blue,thick,ticks=none] table[x=t,y=n,col sep=comma] {beam_periodic_1000_1000_100.dat} ;
\addplot[color=red,thick,dashed,ticks=none] table[x=t,y=e,col sep=comma] {beam_periodic_1000_1000_100.dat} ;
\end{axis}
\end{tikzpicture}
\end{tabular}
\end{center}
\caption{Norm of the solution for a periodic solution with $f=100$Hz for the beam model}
\label{fig:beam-periodic}
\end{figure}

\paragraph{Harmonic solution}
Figure~\ref{fig:beam-periodic} shows the norm of the state vector in function of time for
periodic right-hand side $\mbox{Re}(b_0\cdot e^{\imath2\pi f t})$ with $b_0$ chosen as before.
The starting vector was chosen as in Lemma~\ref{le:case1} for $f=100$Hz.
We chose $n_Z=1000$ sample points for AAA and used three values of $f_{\max}$.
The step size for Crank Nicholson is $1\,10^{-5}$s.
The figure shows the norm of the state vector in a blue solid line and the norm of the absolute error
of the solution in a red dashed line.
We notice that the choice of $f_{\max}$ determines the quality of higher frequency solutions.
This illustrates that $f_{\max}$ should be taken large enough.

\paragraph{Nonzero initial values}
In the following experiment, we used initial values
\[
  \begin{pmatrix}
  x_0 \\ 0 \\ \mathbf{\Phi}(0) \otimes x_0
  \end{pmatrix}.
\]
The vector $x_0 = A(0)^{-1} b_0$, with $b_0$ as defined before,
which corresponds to a constant $x(t)$ for $t\leq 0$.
%We chose $x_0=b_0$ with $b_0$ defined before.
The right-hand side $b(t)$ was chosen identically zero.
As a result, the state vector will gradually decrease to zero for increasing time.
We ran the Crank-Nicholson method with step size $10^{-5}$.
Figure~\ref{fig:beam-constant} shows the norm of $x(t)$ as a function of time obtained for different values of $f_{\max}$ for the AAA approximation.
The tolerance for AAA was set to $10^{-13}$ as before.
We can conclude that it is important to select $f_{\max}$ high enough for obtaining an accurate solution.
Note, however, that the mesh is not valid for such high frequencies, so, we only illustrate the importance of high $f_{\max}$ in the case of sufficiently fine meshes.
\begin{figure}
\begin{center}
\begin{tikzpicture}[scale=1.3]
\begin{axis}[axis lines=center, xlabel={$t$}, ylabel={$\|x(t)\|_2$}]
\addlegendentry{$f_{\max}=10$Hz}
\addplot[color=blue,ticks=none] table[x=t,y=n,col sep=comma] {beam_constant_10.dat} ;
\addlegendentry{$f_{\max}=100$Hz}
\addplot[color=red,ticks=none] table[x=t,y=n,col sep=comma] {beam_constant_100.dat} ;
\addlegendentry{$f_{\max}=10^3$Hz}
\addplot[color=green,dashed,ticks=none] table[x=t,y=n,col sep=comma] {beam_constant_1000.dat} ;
\addlegendentry{$f_{\max}=10^4$Hz}
\addplot[color=brown,dotted,ticks=none] table[x=t,y=n,col sep=comma] {beam_constant_10000.dat} ;
\end{axis}
\end{tikzpicture}
\begin{tikzpicture}[scale=1.3]
\begin{axis}[axis lines=center, xlabel={$t$}, ylabel={$\|x(t)\|_2$}]
\addlegendentry{$f_{\max}=10^3$Hz}
\addplot[color=blue,ticks=none] table[x=t,y=n,col sep=comma] {beam_constant_1000.dat} ;
\addlegendentry{$f_{\max}=10^4$Hz}
\addplot[color=green,dashed,ticks=none] table[x=t,y=n,col sep=comma] {beam_constant_10000.dat} ;
\addlegendentry{$f_{\max}=10^5$Hz}
\addplot[color=red,dashdotted,ticks=none] table[x=t,y=n,col sep=comma] {beam_constant_100000.dat} ;
\end{axis}
\end{tikzpicture}
\end{center}
\caption{Norm of $x(t)$ for initial values associated with constant $x$ for negative times.}\label{fig:beam-constant}
\end{figure}

\begin{figure}
\begin{center}
%\begin{tikzpicture}[scale=1.3]
%\begin{axis}[axis lines=center, xlabel={$t$}, ylabel={$\|x(t)\|_2$}]
%\addlegendentry{$f_{\max}=10$Hz}
%\addplot[color=blue,ticks=none] table[x=t,y=n,col sep=comma] {beam_periodic_initial_10.dat} ;
%\addlegendentry{$f_{\max}=100$Hz}
%\addplot[color=red,ticks=none] table[x=t,y=n,col sep=comma] {beam_periodic_initial_100.dat} ;
%\addlegendentry{$f_{\max}=10^3$Hz}
%\addplot[color=green,dashed,ticks=none] table[x=t,y=n,col sep=comma] {beam_periodic_initial_1000.dat} ;
%\addlegendentry{$f_{\max}=10^4$Hz}
%\addplot[color=brown,dotted,ticks=none] table[x=t,y=n,col sep=comma] {beam_periodic_initial_10000.dat} ;
%\end{axis}
%\end{tikzpicture}
\begin{tikzpicture}[scale=1.3]
\begin{axis}[axis lines=center, xlabel={$t$}, ylabel={$\|x(t)\|_2$}]
\addlegendentry{$f_{\max}=100$Hz}
\addplot[color=blue,ticks=none] table[x=t,y=n,col sep=comma] {beam_periodic_initial_100.dat} ;
\addlegendentry{$f_{\max}=1000$Hz}
\addplot[color=blue,ticks=none] table[x=t,y=n,col sep=comma] {beam_periodic_initial_1000.dat} ;
\addlegendentry{$f_{\max}=10000$Hz}
\addplot[color=orange,ticks=none] table[x=t,y=n,col sep=comma] {beam_periodic_initial_10000.dat} ;
\addlegendentry{$f_{\max}=10^5$Hz}
\addplot[color=red,ticks=none] table[x=t,y=n,col sep=comma] {beam_periodic_initial_100000.dat} ;
\end{axis}
\end{tikzpicture}
\end{center}
\caption{Norm of $x(t)$ for initial values associated with harmonic $x$ for negative times.}\label{fig:beam-periodic_initial}
\end{figure}
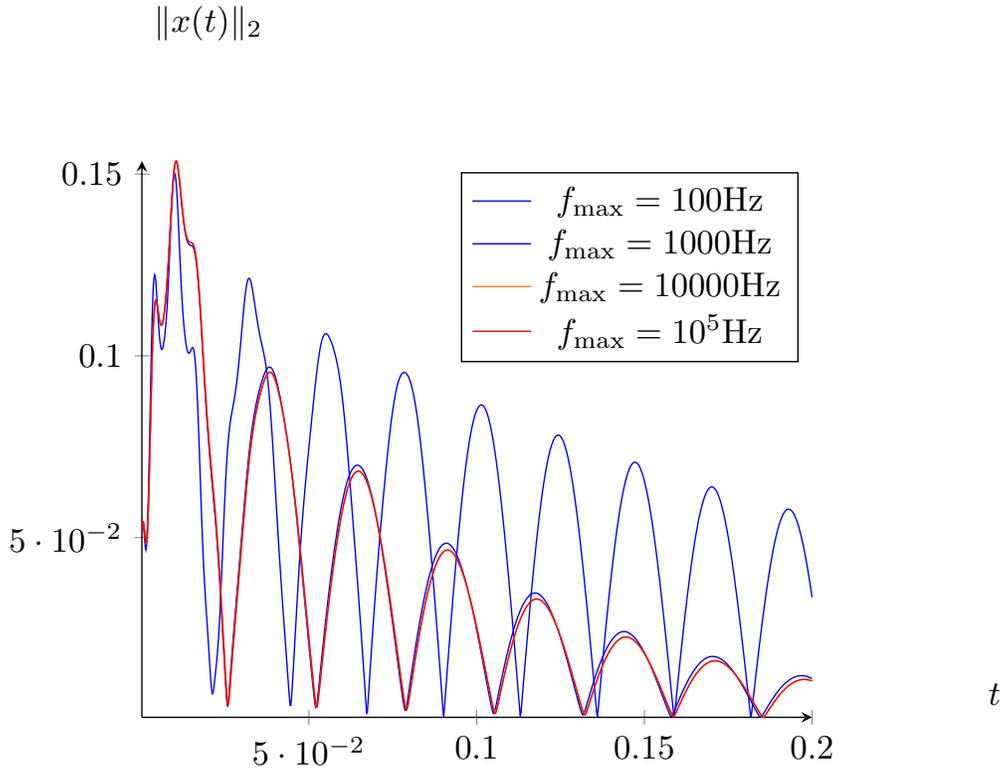

We also ran an experiment with initial values
\[
  \mbox{Re}\begin{pmatrix}
  x_0 \\ \imath\omega x_0 \\ \mathbf{\Phi}(\imath\omega) \otimes x_0
  \end{pmatrix} 
\]
with $\omega=2\pi\cdot100$ which corresponds to a periodic $x(t)$ for $t\leq 0$ with frequency $100$Hz.
The vector $x_0 = A(\i\omega)^{-1} b_0$, with $b_0$ as defined before. 
The right-hand side $b(t)$ was chosen identically zero.
As a result, the state vector will gradually decrease to zero for increasing time.
We ran the Crank-Nicholson method with step size $10^{-5}$.
Figure~\ref{fig:beam-periodic_initial} shows results for three values of $f_{\max}$.
The conclusion is also here that $f_{\max}$ should be chosen large enough for good accuracy.

\paragraph{Stop criteria for AAA}

We now discuss the choice of $\tau$.
We used a harmonic right-hand side and initial solution with $\omega=100$, $f_{\max}=10^4$, $n_Z=10000$ and the Crank-Ncholson method timestep $10^{-5}$.
Since $\|A_{-1}\|\approx 4\,10^2$, $\|A_0\|\approx 2\,10^9$, $\|A_2\|=1\,10^{-3}$ and
$\|g_1\|_Z=4.5\,10^7$, $A_{-1} g_1$ is a dominant term in $A(s)$.

\begin{table}
\begin{center}
\caption{Number of terms, $d$, of the rational approximation for $f_{\max}=10^4$ and $n_Z=10000$ and for different values of AAA's tolerance $\tau$}
\label{tab:number of terms}
\begin{tabular}{cc}\hline
$\tau$ & $d$ \\\hline
$10^{-13}$ & 40 \\
$10^{-10}$ & 28 \\
$10^{-7}$ & 20 \\
$10^{-4}$ & 12 \\
$10^{-1}$ & 4 \\\hline
\end{tabular}
\end{center}
\end{table}

\begin{figure}
\begin{center}
\begin{tabular}{cc}
$\tau=10^{-1}$ & $\tau=10^{-4}$ \\
\begin{tikzpicture}[scale=0.8]
\begin{axis}[axis lines=center, xlabel={$t$}, ylabel={$\|x(t)\|_2$}]
\addlegendentry{Solution}
\addplot[color=blue,thick,ticks=none] table[x=t,y=n,col sep=comma] {beam_tol_1.dat} ;
\addlegendentry{Error Weighted}
\addplot[color=green,thick,dashdotted,ticks=none] table[x=t,y=e,col sep=comma] {beam_tol_1_n.dat} ;
\addlegendentry{Error Set Valued}
\addplot[color=red,thick,dashed,ticks=none] table[x=t,y=e,col sep=comma] {beam_tol_1.dat} ;
\end{axis}
\end{tikzpicture}
&
\begin{tikzpicture}[scale=0.8]
\begin{axis}[axis lines=center, xlabel={$t$}, ylabel={$\|x(t)\|_2$}]
\addlegendentry{Solution}
\addplot[color=blue,thick,ticks=none] table[x=t,y=n,col sep=comma] {beam_tol_4.dat} ;
\addlegendentry{Error Weighted}
\addplot[color=green,thick,dashdotted,ticks=none] table[x=t,y=e,col sep=comma] {beam_tol_4_n.dat} ;
\addlegendentry{Error Set Valued}
\addplot[color=red,thick,dashed,ticks=none] table[x=t,y=e,col sep=comma] {beam_tol_4.dat} ;
\end{axis}
\end{tikzpicture}
\end{tabular}
\end{center}
\caption{Norm of the solution for a harmonic solution with $f=100$Hz}
\label{fig:beam-periodic-tol}
\end{figure}
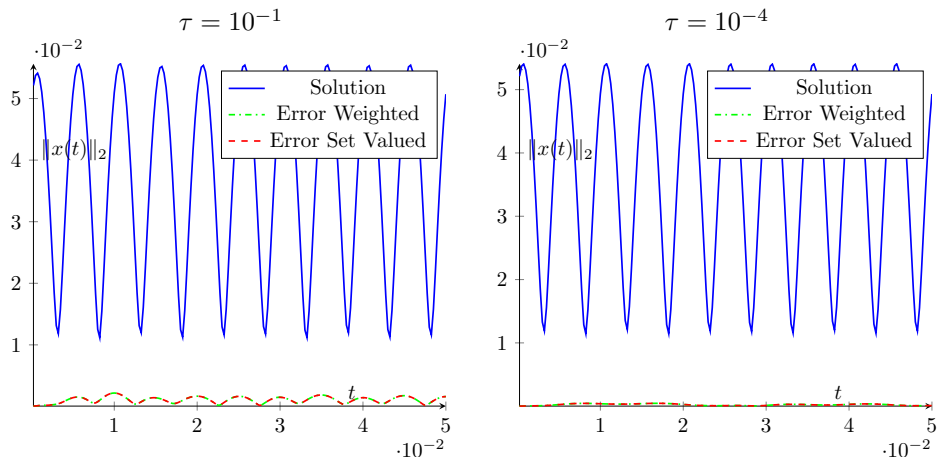

Table~\ref{tab:number of terms} shows the number of terms in the rational approximation as a function of the tolerance $\tau$.
Figure~\ref{fig:beam-periodic-tol} shows the error of the solution for two choices of $\tau$.
For $\tau=10^{-4}$, we have a solution of excellent quality, with $d=12$.
Even for $d=10^{-1}$, the quality of the solution is acceptable.

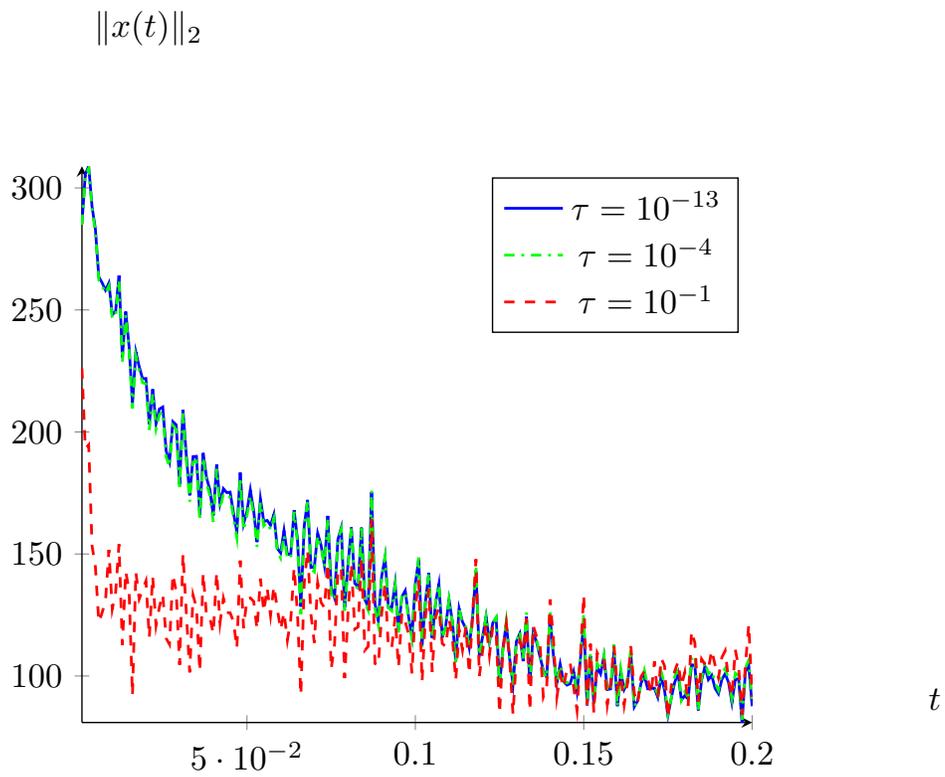
\begin{figure}
\begin{center}
\begin{tikzpicture}[scale=1.3]
\begin{axis}[axis lines=center, xlabel={$t$}, ylabel={$\|x(t)\|_2$}]
\addlegendentry{$\tau=10^{-13}$}
\addplot[color=blue,thick,ticks=none] table[x=t,y=n,col sep=comma] {beam_constant_tol_13.dat} ;
\addlegendentry{$\tau=10^{-4}$}
\addplot[color=green,thick,dashdotted,ticks=none] table[x=t,y=n,col sep=comma] {beam_constant_tol_4.dat} ;
\addlegendentry{$\tau=10^{-1}$}
\addplot[color=red,thick,dashed,ticks=none] table[x=t,y=n,col sep=comma] {beam_constant_tol_1.dat} ;
\end{axis}
\end{tikzpicture}
\end{center}
\caption{Norm of the solution for a zero right hand side}
\label{fig:beam-zero-tol}
\end{figure}

\subsection{Porous car seats}
\begin{figure}[!h]
\begin{center}
%\centering
\begin{subfigure}{0.45\textwidth}
\includegraphics[width=\textwidth]{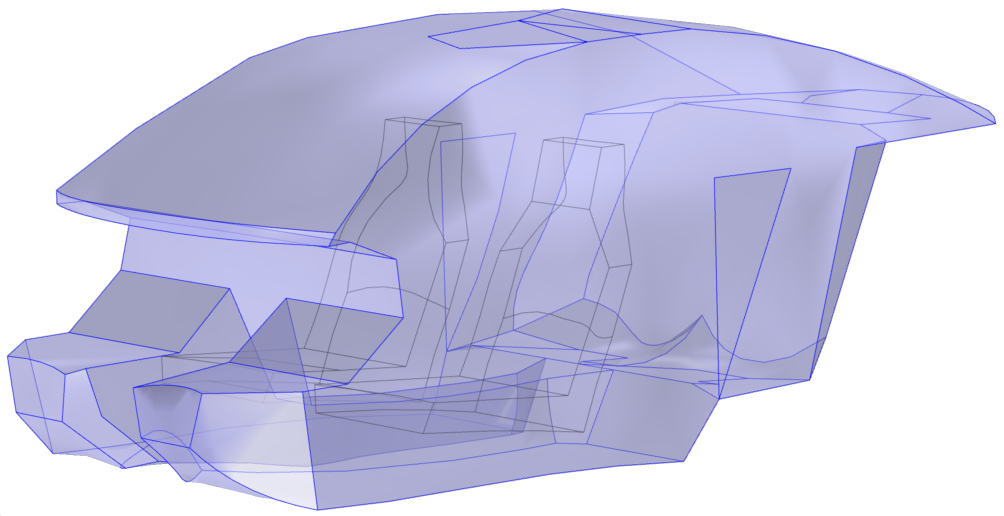}
   \label{fig:carcavitygeom}
\end{subfigure}
%\centering
\begin{subfigure}{0.45\textwidth}
\includegraphics[width=\textwidth]{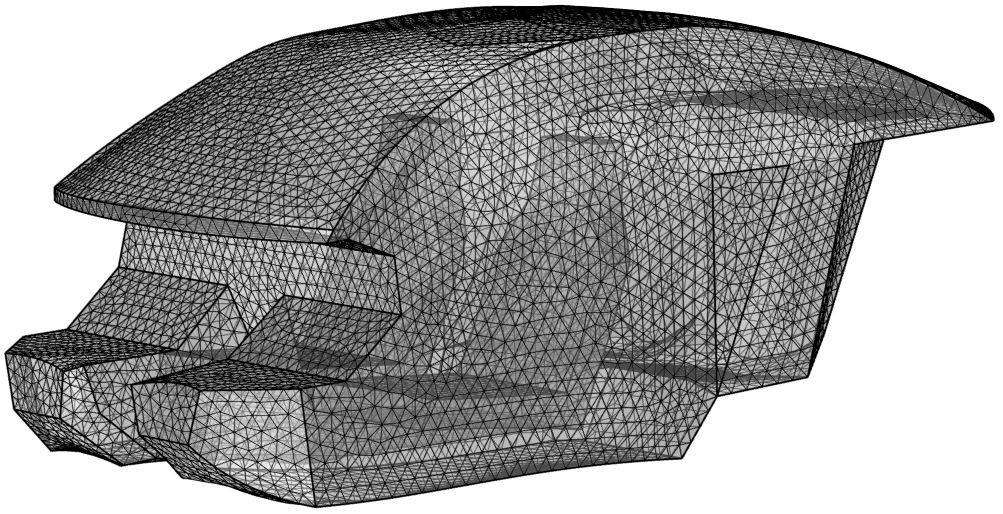}
   \label{fig:carcavitymesh}
\end{subfigure}
\end{center}
\caption{Second application case: Set-up of the poro-acoustic car interior geometry. Left geometry, right FE mesh}
\label{fig:carcavity}
\end{figure}

The second problem case considers a porous-acoustic problem, consisting of an acoustic car interior geometry, with two seats as shown in figure \ref{fig:carcavity}. The air inside the cavity is modelled using the acoustic Helmholtz equation assuming a fluid density of $\rho_0=1.213kg/m^3$ and a speed of sound $c_0 = 342.0m/s$. The Johnson-Champoux-Allard rigid frame equivalent fluid model \cite{Allard2009} is used to describe the frequency-dependent behaviour of the porous seats, accounting for thermal and viscous losses.\\ 

Starting from the acoustic mass, $M_a$ and $M_p$, and stiffness matrices, $K_a$ and $K_p$, with subscript $a$ and $p$ representing the acoustic and porous domain, the complex system matrix $A(s)$, is constructed as follows: 
\[
  A(s) = (K_a + K_p g_1(s)/\phi)/\rho_0 + s^2 (M_a + M_p g_2(s)/\phi)/ (\rho_0 c_0^2)
\]
with
\begin{eqnarray*}
g_1(s) & = & \phi / \alpha(s) \\
\alpha(s) & = & \alpha_{\infty} (1+(\sigma\phi/(s \rho_0 \alpha_{\infty}))G_J(s)) \\
G_J(s) & = & \sqrt{1+(4 \alpha_{\infty}^2\eta\rho_0s)/(\sigma^2\Lambda^2\phi^2)}
\end{eqnarray*}
and
\begin{eqnarray*}
g_2 & = & \phi (\gamma - (\gamma - 1)/ \alpha'(s)) \\
\alpha'(s) & = & 1 + 8\eta/(\Lambda'^2 N_{\rm Pr} s\rho_0) \sqrt{1 + \rho_0 s N_{\rm Pr}\Lambda'^2/(16 \eta)}.
\end{eqnarray*}
Table~\ref{tab:constants car} gives the meaning of the parameters and the concrete values used in this specific example.
\begin{table}
%\begin{center}
%\begin{tabular}{|c|c|c|}\hline
%$\phi = 0.98$  & $\rho_0 = 1.213$           & $P_0 = 101325$ \\\hline
%$\gamma = 1.4$   & $\eta = 1.839\,10^{-5}$ & $R_{\mbox{gas}} = 287.031$  \\\hline
%$T_{\rm absolute} = P_0/(\rho_0 R_{\rm gas})$ & $c_0 = \sqrt{\gamma R_{\rm gas} T_{\rm absolute}}$ & $\sigma = 13500$ \\\hline
%\multicolumn{3}{|c|}{$C_p = 4168.8 0.249679-7.55179\,10^{-5}  T_{\rm absolute}+1.69194\,10^{-7}T_{\rm absolute}^2-6.46128\,10^{-11} T_{\rm absolute}^3$} \\\hline
%\multicolumn{3}{|c|}{$\kappa = 2.624\,10^{-2} ((T_{\rm absolute}/300)^{1.5} (300+245.4\exp(-27.6/300)))/(T_{\rm absolute}+245.4\exp(-27.6/T_{\rm absolute}))$}
%\\\hline
%$\phi = 0.98$ & $\Lambda = 80\,10^{-6}$ & $\Lambda' = 160\,10^{-6}$ \\\hline
%$\alpha_{\infty} = 1.7$ & $N_{\rm Pr} = \eta C_p/\kappa$ & \\\hline
%\end{tabular}
%\end{center}
%\caption{Constants for the porous car seat problem}\label{tab:constants car}
  \centering
  \begin{tabular}{ccc|ccc} \toprule
    \multicolumn{6}{c}{\textbf{Material properties}} \\ \cmidrule(r){1-6}  
    \multicolumn{3}{c|}{\textbf{Properties of air}} & 
    \multicolumn{3}{c}{\textbf{Properties of porous material}} \\
    \cmidrule(r){1-6}
    $\rho_0$ &  air density  & $1.213$ kg/m$^3$ & $\phi$ & porosity & $0.98$ \\ 
    $\mathrm{Pr}$ &  Prandtl number  & $0.72$ & $\alpha_{\infty}$ & tortuosity & $1.7$\\
    $\gamma$ & heat capacity ratio  & $1.4$ & $\sigma$ & flow resistivity & $13500$ Ns/m$^4$\\
    $\eta$ &  dynamic viscosity  & $0.1837$ $\mu$kg/(m$\cdot$s) & $\Lambda$ & viscous length & $80$ $\mu$m\\
    & & & $\Lambda'$ & thermal length & $160$ $\mu$m\\
    \midrule
  \end{tabular}
  \caption{{Material properties of the air and the porous material applied in the acoustic car interior geometry with porous seats.}}
\label{tab:constants car}
\end{table}
The applied damping model is known to be causal, a time domain formulation has been presented in \cite{umnova2009}.
The complex function $\alpha(s)$ is called the dynamic tortuosity and is used to assess the effective density in the porous material via $\rho(s)=\alpha(s)\rho_0$.
The expression was derived by Johnson et al.\ \cite{Johnson1987} who justified its use by physical causality constraints concerning its singularities
which must be located on the negative real axis \cite{Allard2009} and correct asymptotic behaviour for low and high frequencies:
from microscale Stokes flow at very low frequencies to inviscid flow as high frequency asymptote.
Function $g_1$ indeed has a pole at $-\sigma\phi/\rho_0\alpha_{\infty}$ and a singularity at $-\sigma^2\Lambda^2\phi^2/4\alpha_{\infty}^2\eta\rho_0$, which both lie on the negative real axis.
Analogously, the complex function $\alpha'(s)$ was developed by Lafarge et al.\ \cite{lafarge1997} to arrive at an equivalent bulk density $K=P_0/(1-\frac{\gamma-1}{\gamma\alpha'(s)})$,
which accounts for the transition from isothermal behaviour at low frequencies to adiabatic behaviour as high-frequency assymptote.
Function $g_2$ has a singularity at $-16\eta/\rho_0N_{\rm Pr}\Lambda'^2$, and there is a pole at $s=-(\sqrt{17}-1)2\eta/\Lambda'^2 N_{\rm Pr} \rho_0$.

%\begin{table}
%\begin{center}
%\begin{tabular}{|cc|cc|cc|}\hline
%$\phi$   & $0.98$  & $\rho_0$ & $1.213$           & $P_0$   & $101325$ \\\hline
%$\gamma$ & $1.4$   & $\eta$   &  $1.839\,10^{-5}$ & $rho_0$ & $1.213$ \\\hline
%$R_{\mbox{gas}}$ & $287.031$ &  $T_{\rm absolute}$ & $P_0/(rho_0 R_{\rm gas})$ & $c_0$ %& $\sqrt{\gamma R_{\rm gas} T_{\rm absolute}}$ \\\hline
%$C_p$ & \multicolumn{5}{c|}{$4168.8 0.249679-7.55179\,10^{-5}  T_{\rm %absolute}+1.69194\,10^{-7}T_{\rm absolute}^2-6.46128\,10^{-11} T_{\rm absolute}^3$} %\\\hline
%$\kappa$ & \multicolumn{3}{c|}{$2.624\,10^{-2} \frac{ (T_{\rm absolute}/300)^{1.5} %(300+245.4\exp(-27.6/300))}{T_{\rm absolute}+245.4\exp(-27.6/T_{\rm absolute})}$} &
%$N_{\rm Pr}$ &  $\eta C_p/\kappa$ \\\hline
%$\sigma$ & $13500$ & $\phi$ & $0.98$ & $\Lambda$ & $80\,10^{-6}$ \\\hline
%$\Lambda'$ & $160\,10^{-6}$ & $\alpha_{\infty}$ & $1.7$ & & \\\hline
%\end{tabular}
%\end{center}
%\caption{Constants for the porous car seat problem}\label{tab:constants car}
%\end{table}
The matrices $K_a$, $K_p$, $M_a$ and $M_p$ are sparse real symmetric matrices of order 81,257.

\paragraph{AAA approximation} We first show the quality of the AAA algorithm and its variants.
%The stability was checked by verifying the signs of the real parts of the poles of the rational approximation.
The weighted AAA method takes into account the relative contribution of the two functions $g_1$ and $g_2$ to the
nonlinear matrix.
We found that the functions $g_1$ and $g_2=s^2 \tilde{g}_2(s)$ do mostly not lead to stable poles for any $f_{\max}=10$, $100$,
$1000$ and $10,000$.
For example, for $f_{\max}=100$, $N_Z=1000$ and $\tau=10^{-12}$, we found the poles
\[
509981,-397518,-19086.3,-8148.97,-5125.9,-3237.91,-1951.31.
\]
For $f_{\max}=1000$, $N_Z=10000$ and $\tau=10^{-12}$, we found that 7 out of 23 poles were unstable.

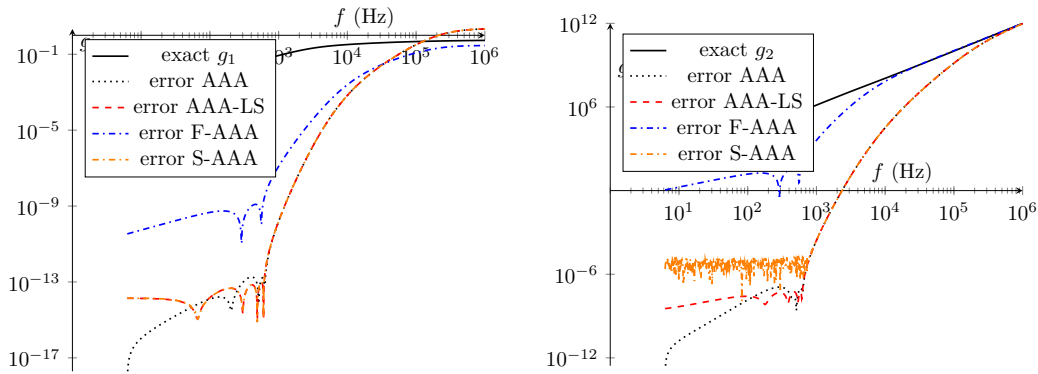
\begin{figure}
\begin{center}
\begin{tabular}{cc}
\begin{tikzpicture}[scale=0.8]
\begin{loglogaxis}[axis lines=center, xlabel={$f$ (Hz)}, ylabel={$g_1$},xmin=1,xmax=1000000,legend pos=north west]
\addlegendentry{exact $g_1$}
\addplot[color=black,thick,ticks=none] table[x=s,y=n1,col sep=comma] {porous2/aaa_diff_100_1000.dat} ;
\addlegendentry{error AAA}
\addplot[color=black,dotted,thick,ticks=none] table[x=s,y=e1,col sep=comma] {porous2/aaa_diff_100_1000.dat} ;
\addlegendentry{error AAA-LS}
\addplot[color=red,thick,dashed,ticks=none] table[x=s,y=e1,col sep=comma] {porous2/aaa_ls_diff_100_1000.dat} ;
\addlegendentry{error F-AAA}
\addplot[color=blue,thick,dashdotted,ticks=none] table[x=s,y=e1,col sep=comma] {porous2/f_aaa_diff_100_1000.dat} ;
\addlegendentry{error S-AAA}
\addplot[color=orange,thick,dashdotted,ticks=none] table[x=s,y=e1,col sep=comma] {porous2/s_aaa_diff_100_1000.dat} ;
\end{loglogaxis}
\end{tikzpicture} &
\begin{tikzpicture}[scale=0.8]
\begin{loglogaxis}[axis lines=center, xlabel={$f$ (Hz)}, ylabel={$g_2$},xmin=1,xmax=1000000,legend pos=north west]
\addlegendentry{exact $g_2$}
\addplot[color=black,thick,ticks=none] table[x=s,y=n2,col sep=comma] {porous2/aaa_diff_100_1000.dat} ;
\addlegendentry{error AAA}
\addplot[color=black,thick,dotted,ticks=none] table[x=s,y=e2,col sep=comma] {porous2/aaa_diff_100_1000.dat} ;
\addlegendentry{error AAA-LS}
\addplot[color=red,thick,dashed,ticks=none] table[x=s,y=e2,col sep=comma] {porous2/aaa_ls_diff_100_1000.dat} ;
\addlegendentry{error F-AAA}
\addplot[color=blue,thick,dashdotted,ticks=none] table[x=s,y=e2,col sep=comma] {porous2/f_aaa_diff_100_1000.dat} ;
\addlegendentry{error S-AAA}
\addplot[color=orange,thick,dashdotted,ticks=none] table[x=s,y=e2,col sep=comma] {porous2/s_aaa_diff_100_1000.dat} ;
\end{loglogaxis}
\end{tikzpicture}
\end{tabular}
\end{center}
\caption{Comparison of AAA, AAA-LS, F-AAA, and S-AAA for $g_1$ and $g_2$ for $f_{\max}=100$ and $N_Z=1000$.}\label{fig:porous-aaa-ls}
\end{figure}

With the filtered versions of AAA, the unstable poles are eliminated in order to obtain a stable approximation.
We used both the representation using weighted partial fractions and weighted inverted Newton polynomials.
We did not observe a difference in accuracy.
Therefore, we only report results with weighted partial fractions.
The elimination of the unstable poles and the reduction of the number of rational terms
reduces the accuracy of the approximation as can be seen in Figure~\ref{fig:porous-aaa-ls}.
For this problem, flipping the unstable pole to the left half plane shows good accuracy.

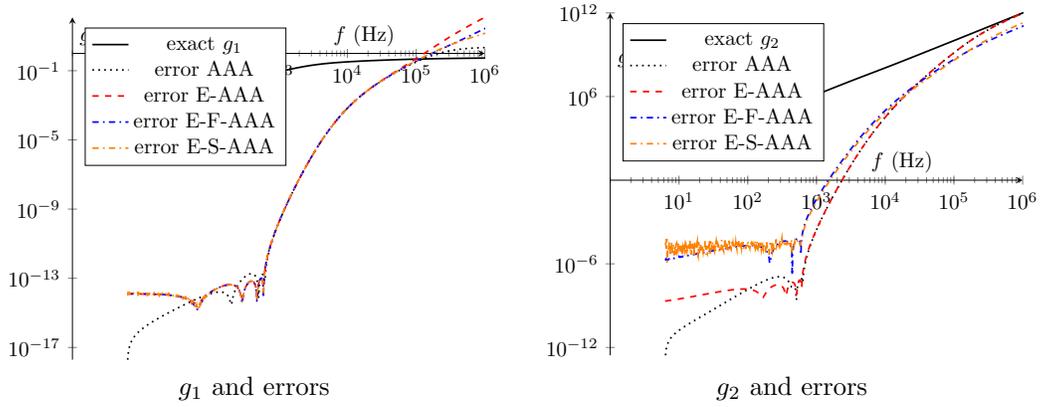
\begin{figure}
\begin{center}
\begin{tabular}{cc}
\begin{tikzpicture}[scale=0.8]
\begin{loglogaxis}[axis lines=center, xlabel={$f$ (Hz)}, ylabel={$g_1$},xmin=1,xmax=1000000,legend pos=north west]
\addlegendentry{exact $g_1$}
\addplot[color=black,thick,ticks=none] table[x=s,y=n1,col sep=comma] {porous2/aaa_diff_100_1000.dat} ;
\addlegendentry{error AAA}
\addplot[color=black,dotted,thick,ticks=none] table[x=s,y=e1,col sep=comma] {porous2/aaa_diff_100_1000.dat} ;
\addlegendentry{error E-AAA}
\addplot[color=red,thick,dashed,ticks=none] table[x=s,y=e1,col sep=comma] {porous2/e_aaa_diff_100_1000.dat} ;
\addlegendentry{error E-F-AAA}
\addplot[color=blue,thick,dashdotted,ticks=none] table[x=s,y=e1,col sep=comma] {porous2/e_f_aaa_diff_100_1000.dat} ;
\addlegendentry{error E-S-AAA}
\addplot[color=orange,thick,dashdotted,ticks=none] table[x=s,y=e1,col sep=comma] {porous2/e_s_aaa_diff_100_1000.dat} ;
\end{loglogaxis}
\end{tikzpicture} &
\begin{tikzpicture}[scale=0.8]
\begin{loglogaxis}[axis lines=center, xlabel={$f$ (Hz)}, ylabel={$g_2$},xmin=1,xmax=1000000,legend pos=north west]
\addlegendentry{exact $g_2$}
\addplot[color=black,thick,ticks=none] table[x=s,y=n2,col sep=comma] {porous2/aaa_diff_100_1000.dat} ;
\addlegendentry{error AAA}
\addplot[color=black,thick,dotted,ticks=none] table[x=s,y=e2,col sep=comma] {porous2/aaa_diff_100_1000.dat} ;
\addlegendentry{error E-AAA}
\addplot[color=red,thick,dashed,ticks=none] table[x=s,y=e2,col sep=comma] {porous2/e_aaa_diff_100_1000.dat} ;
\addlegendentry{error E-F-AAA}
\addplot[color=blue,thick,dashdotted,ticks=none] table[x=s,y=e2,col sep=comma] {porous2/e_f_aaa_diff_100_1000.dat} ;
\addlegendentry{error E-S-AAA}
\addplot[color=orange,thick,dashdotted,ticks=none] table[x=s,y=e2,col sep=comma] {porous2/e_s_aaa_diff_100_1000.dat} ;
\end{loglogaxis}
\end{tikzpicture} \\
    $g_1$ and errors & $g_2$ and errors
\end{tabular}
\end{center}
\caption{Comparison of AAA, E-AAA, E-S-AAA, and F-E-AAA for $g_1$ and $g_2$ for $f_{\max}=100$Hz and $N_Z=1000$.}\label{fig:porous-e-aaa-ls}
\end{figure}

Figure~\ref{fig:porous-aaa-ls} shows results for weighted AAA, AAA-LS, S-AAA, and F-AAA.
Since weighted AAA produces unstable poles, the filter step is advisable for time integration.
We also compare with extended AAA.
This approach appears to suffer much less from the reduction in accuracy as is shown in Figure~\ref{fig:porous-e-aaa-ls}.

We also tried other choices of functions to build the AAA approximation.
Since the extended AAA allows us to make an approximation of $s g$ and $s^2 g$ from a AAA approximation of $g$, we show results for
applying AAA to $[g_1,g_2/s]$ and $[g_1,g_2/s^2]$.

\begin{figure}
\begin{center}
\begin{tabular}{cc}
\begin{tikzpicture}[scale=0.8]
\begin{loglogaxis}[axis lines=center, xlabel={$f$ (Hz)}, ylabel={$g_1$},xmin=1,xmax=1000000,legend pos=north west]
\addlegendentry{exact $g_1$}
\addplot[color=black,thick,ticks=none] table[x=s,y=n1,col sep=comma] {porous2/aaa_diff_100_1000.dat} ;
\addlegendentry{error AAA}
\addplot[color=black,thick,dotted,ticks=none] table[x=s,y=e1,col sep=comma] {porous2/aaa_diff_100_1000_1.dat} ;
\addlegendentry{error AAA-LS}
\addplot[color=black,thick,dashed,ticks=none] table[x=s,y=e1,col sep=comma] {porous2/aaa_ls_diff_100_1000_1.dat} ;
\addlegendentry{error E-AAA}
\addplot[color=red,thick,dashed,ticks=none] table[x=s,y=e1,col sep=comma] {porous2/e_aaa_diff_100_1000_1.dat} ;
\addlegendentry{error E-F-AAA}
\addplot[color=blue,thick,dashdotted,ticks=none] table[x=s,y=e1,col sep=comma] {porous2/e_f_aaa_diff_100_1000_1.dat} ;
\end{loglogaxis}
\end{tikzpicture} &
\begin{tikzpicture}[scale=0.8]
\begin{loglogaxis}[axis lines=center, xlabel={$f$ (Hz)}, ylabel={$g_2$},xmin=1,xmax=1000000,legend pos=north west]
\addlegendentry{exact $g_2$}
\addplot[color=black,thick,ticks=none] table[x=s,y=n2,col sep=comma] {porous2/aaa_diff_100_1000.dat} ;
\addlegendentry{error AAA}
\addplot[color=black,thick,dotted,ticks=none] table[x=s,y=e2,col sep=comma] {porous2/aaa_diff_100_1000_1.dat} ;
\addlegendentry{error AAA-LS}
\addplot[color=black,thick,dashed,ticks=none] table[x=s,y=e2,col sep=comma] {porous2/aaa_ls_diff_100_1000_1.dat} ;
\addlegendentry{error E-AAA}
\addplot[color=red,thick,dashed,ticks=none] table[x=s,y=e2,col sep=comma] {porous2/e_aaa_diff_100_1000_1.dat} ;
\addlegendentry{error E-F-AAA}
\addplot[color=blue,thick,dashdotted,ticks=none] table[x=s,y=e2,col sep=comma] {porous2/e_f_aaa_diff_100_1000_1.dat} ;
\end{loglogaxis}
\end{tikzpicture} \\
\end{tabular}
\end{center}
\caption{Comparison of AAA, AAA-LS and filtered AAA-LS for AAA on $[g_1,g_2/s]$ for $f_{\max}=100$Hz and $N_Z=1000$.}\label{fig:porous-e1-aaa-ls}
\end{figure}
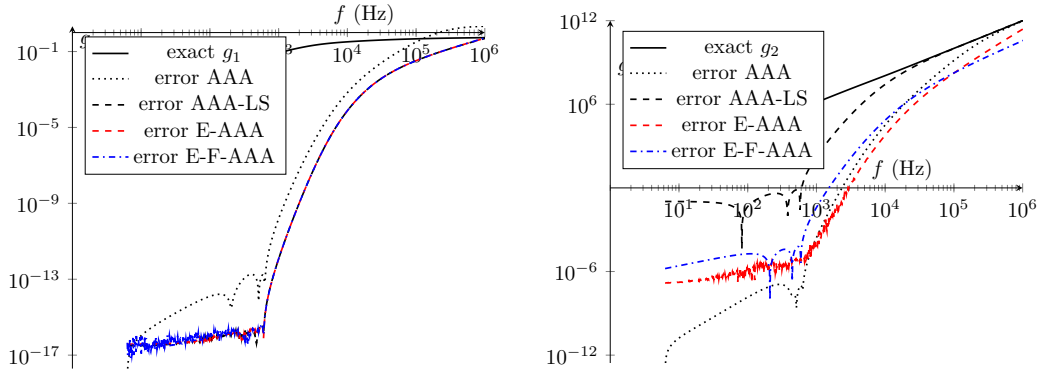

The function $g_1$ is well approximated.
Function $g_2$ contains a factor $s^2$, which has a similar effect as a mass term.
Therefore, we decided to use a factor $s$ instead to make the function look more like a damping term.
That is, we apply weighted AAA to the set $\{g_1,g_2/s\}$ instead of $\{g_1,g_2\}$.
Figure~\ref{fig:porous-e1-aaa-ls} shows the difference in approximation between $g_2$ and $g_2/s$.
The line for AAA is the result for AAA applied to $\{g_1,g_2\}$. It is added for comparisons.
It does not lead to a better approximation.
We also see that the extension of AAA by a polynomial is needed to capture $g_2$ well: AAA-LS shows a better error for $g_1$ compared to AAA, but not for $g_2$.
We obtained more often stable approximations of $\{g_1,g_2/s\}$ than of $\{g_1,g_2\}$.

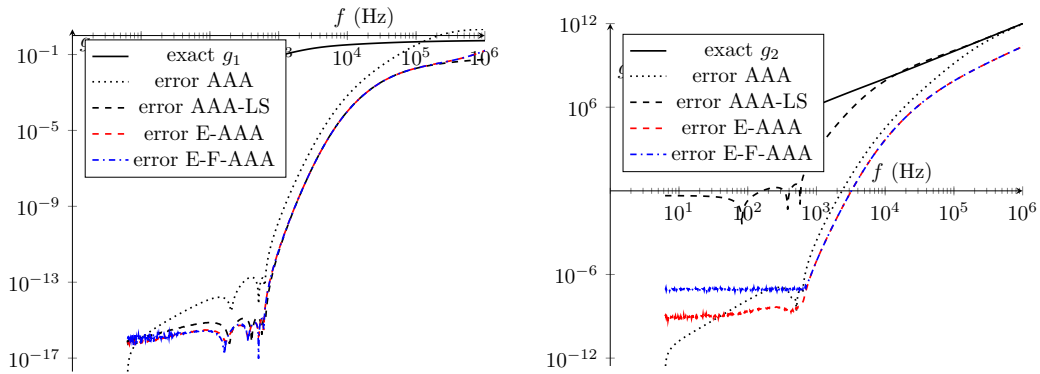
\begin{figure}
\begin{center}
\begin{tabular}{cc}
\begin{tikzpicture}[scale=0.8]
\begin{loglogaxis}[axis lines=center, xlabel={$f$ (Hz)}, ylabel={$g_1$},xmin=1,xmax=1000000,legend pos=north west]
\addlegendentry{exact $g_1$}
\addplot[color=black,thick,ticks=none] table[x=s,y=n1,col sep=comma] {porous2/aaa_diff_100_1000.dat} ;
\addlegendentry{error AAA}
\addplot[color=black,thick,dotted,ticks=none] table[x=s,y=e1,col sep=comma] {porous2/aaa_diff_100_1000_2.dat} ;
\addlegendentry{error AAA-LS}
\addplot[color=black,thick,dashed,ticks=none] table[x=s,y=e1,col sep=comma] {porous2/aaa_ls_diff_100_1000_2.dat} ;
\addlegendentry{error E-AAA}
\addplot[color=red,thick,dashed,ticks=none] table[x=s,y=e1,col sep=comma] {porous2/e_aaa_diff_100_1000_2.dat} ;
\addlegendentry{error E-F-AAA}
\addplot[color=blue,thick,dashdotted,ticks=none] table[x=s,y=e1,col sep=comma] {porous2/e_f_aaa_diff_100_1000_2.dat} ;
\end{loglogaxis}
\end{tikzpicture} &
\begin{tikzpicture}[scale=0.8]
\begin{loglogaxis}[axis lines=center, xlabel={$f$ (Hz)}, ylabel={$g_2$},xmin=1,xmax=1000000,legend pos=north west]
\addlegendentry{exact $g_2$}
\addplot[color=black,thick,ticks=none] table[x=s,y=n2,col sep=comma] {porous2/aaa_diff_100_1000.dat} ;
\addlegendentry{error AAA}
\addplot[color=black,thick,dotted,ticks=none] table[x=s,y=e2,col sep=comma] {porous2/aaa_diff_100_1000_2.dat} ;
\addlegendentry{error AAA-LS}
\addplot[color=black,thick,dashed,ticks=none] table[x=s,y=e2,col sep=comma] {porous2/aaa_ls_diff_100_1000_2.dat} ;
\addlegendentry{error E-AAA}
\addplot[color=red,thick,dashed,ticks=none] table[x=s,y=e2,col sep=comma] {porous2/e_aaa_diff_100_1000_2.dat} ;
\addlegendentry{error E-F-AAA}
\addplot[color=blue,thick,dashdotted,ticks=none] table[x=s,y=e2,col sep=comma] {porous2/e_f_aaa_diff_100_1000_2.dat} ;
\end{loglogaxis}
\end{tikzpicture} \\
\end{tabular}
\end{center}
\caption{Comparison of AAA, AAA-LS and filtered AAA-LS for AAA on $[g_1,g_2/s^2]$ for $f_{\max}=100$Hz and $N_Z=1000$.}\label{fig:porous-e2-aaa-ls}
\end{figure}

\begin{figure}
\begin{center}
\begin{tabular}{cc}
\begin{tikzpicture}[scale=0.8]
\begin{loglogaxis}[axis lines=center, xlabel={$f$ (Hz)}, ylabel={$g_1$},xmin=1,xmax=1000000,legend pos=north west]
\addlegendentry{exact $g_2$}
\addplot[color=black,thick,ticks=none] table[x=s,y=n2,col sep=comma] {porous2/aaa_diff_100_1000.dat} ;
\addlegendentry{error AAA}
\addplot[color=black,thick,dotted,ticks=none] table[x=s,y=e2,col sep=comma] {porous2/aaa_diff_100_1000_2.dat} ;
\addlegendentry{error E-AAA}
\addplot[color=red,thick,dashed,ticks=none] table[x=s,y=e2,col sep=comma] {porous2/e_aaa_diff_100_1000_2.dat} ;
\addlegendentry{error E-F-AAA}
\addplot[color=blue,thick,dashdotted,ticks=none] table[x=s,y=e2,col sep=comma] {porous2/e_f_aaa_diff_100_1000_2.dat} ;
\addlegendentry{error E-S-AAA}
\addplot[color=orange,thick,dashdotted,ticks=none] table[x=s,y=e2,col sep=comma] {porous2/e_s_aaa_diff_100_1000_2.dat} ;
\end{loglogaxis}
\end{tikzpicture} &
\begin{tikzpicture}[scale=0.8]
\begin{loglogaxis}[axis lines=center, xlabel={$f$ (Hz)}, ylabel={$g_2$},xmin=1,xmax=1000000,legend pos=north west]
\addlegendentry{exact $g_2$}
\addplot[color=black,thick,ticks=none] table[x=s,y=n2,col sep=comma] {porous2/aaa_diff_1000_10000_2.dat} ;
\addlegendentry{error AAA}
\addplot[color=black,thick,dotted,ticks=none] table[x=s,y=e2,col sep=comma] {porous2/aaa_diff_1000_10000_2.dat} ;
\addlegendentry{error E-AAA}
\addplot[color=red,thick,dashed,ticks=none] table[x=s,y=e2,col sep=comma] {porous2/e_aaa_diff_1000_10000_2.dat} ;
\addlegendentry{error E-F-AAA}
\addplot[color=blue,thick,dashdotted,ticks=none] table[x=s,y=e2,col sep=comma] {porous2/e_f_aaa_diff_1000_10000_2.dat} ;
\addlegendentry{error E-S-AAA}
\addplot[color=orange,thick,dashdotted,ticks=none] table[x=s,y=e2,col sep=comma] {porous2/e_s_aaa_diff_1000_10000_2.dat} ;
\end{loglogaxis}
\end{tikzpicture} \\
    $f_{\max}=100$ & $f_{\max}=1000$
\end{tabular}
\end{center}
\caption{Comparison of AAA, AAA-LS and filtered AAA-LS for AAA on $[g_1,g_2/s^2]$ for $f_{\max}=100$Hz, $f_{\max}=1000$Hz and $N_Z=10000$.}\label{fig:porous-aaa-100-1000}
\end{figure}
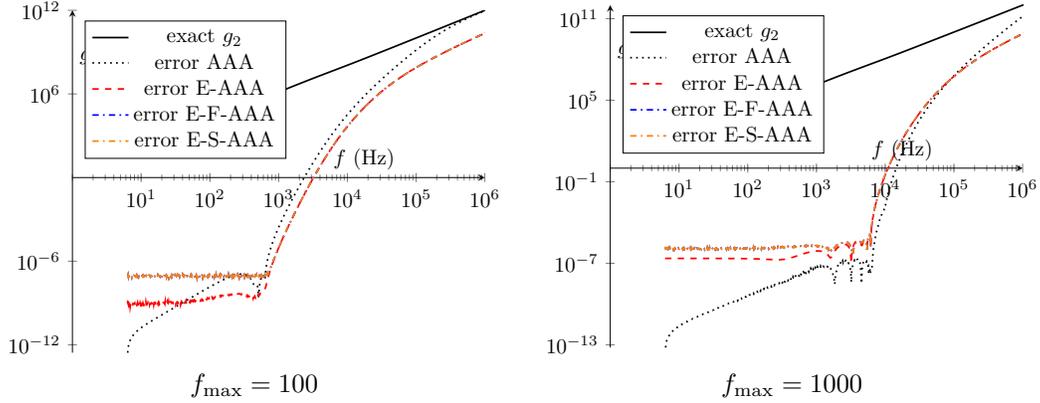

The second choice, $\{g_1,g_2/s^2\}$, is particularly appealing because $g_2/s^2$ is the factor added to the mass matrix, and is a function that varies slowly.
In Figure~\ref{fig:porous-e2-aaa-ls},
also here, AAA-LS does not show good results, where the extension with a polynomial term, as in E-AAA, produces a small error.
the choice $\{g_1,g_2/s^2\}$ shows an improvement of the error for E-F-AAA-LS.
Figure~\ref{fig:porous-aaa-100-1000} shows results of the approximation of $g_2$ for two frequencies, 100Hz and 1000Hz.

%\begin{tabular}{cccc}
%Method & $[g_1,g_2]$ & $[g_1,g_2/s]$ & $[g_1,g_2/s^2]$ \\\hline
%AAA    & 6 & 6 & 6 \\
%E-AAA-LS & 6 6 & 6 \\
%EF-AAA-LS & 4 & 5 & 5
%\end{tabular}

\begin{figure}
\begin{center}
\begin{tikzpicture}[scale=0.8]
\begin{semilogyaxis}[axis lines=center, xlabel={$t$}, ylabel={$\|x(t)\|_2$}, ymax=50, legend pos=south east]
\addlegendentry{Solution}
\addplot[color=black,thick,ticks=none] table[x=t,y=n,col sep=comma] {porous2/filter_periodic_100_10000_100.dat} ;
\addlegendentry{Error $f_{\max}=100$Hz}
\addplot[color=blue,thick,dashed,ticks=none] table[x=t,y=e,col sep=comma] {porous2/filter_periodic_100_10000_100.dat} ;
\addlegendentry{Error $f_{\max}=1000$Hz}
\addplot[color=brown,thick,dashdotted,ticks=none] table[x=t,y=e,col sep=comma] {porous2/filter_periodic_1000_10000_100.dat} ;
\addlegendentry{Error $f_{\max}=10000$Hz}
\addplot[color=brown,thick,dashdotted,ticks=none] table[x=t,y=e,col sep=comma] {porous2/filter_periodic_10000_10000_100.dat} ;
\end{semilogyaxis}
\end{tikzpicture}
\end{center}
\caption{Norm of the solution and the error for a harmonic solution with $f=100$Hz for the porous car seat model using the F-E-AAA-LS approximation method}
\label{fig:porous-periodic}
\end{figure}
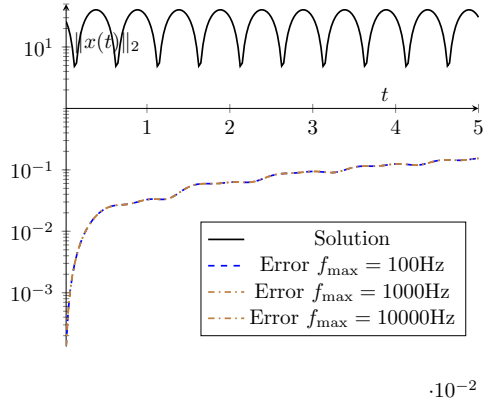

\paragraph{Periodic solution}
We tried several values of $f_{\max}$ for a periodic solution as shown in Lemma~\ref{le:case1} for frequency $100$Hz.
The right-hand side is chosen as $f_0\cos(\omega t)$ with $\omega=2\pi 100$. %%%%%!!!!, and $f_0$ ...
Figure~\ref{fig:porous-periodic} shows the results for various choices of $f_{\max}$, $\tau=10^{-13}$ and $n_Z=10000$.
Standard weighted AAA and the extended version did not produce a stable solution.
Filtering is therefore needed.
We have used F-E-AAA-LS for time integration.
With $f_{\max}=20$Hz, we obtained an approximation with stable poles, but time integration was not stable.
This may not be surprising because for such a low frequency range, the eigenvalues corresponding with high frequencies are
badly approximated and may be unstable.
For $f_{\max}=100$Hz and $f_{\max}=1000$Hz, a stable solution was obtained, as can be seen in Figure~\ref{fig:porous-periodic}.

\section{Conclusions}\label{sec:conc}

We have discussed the development of linear systems whose output approximates the output of a system with nonlinear frequency dependencies, and their use in time integration.
We showed the connection of the nonlinear matrix with the Laplace transform of the linearization and gave numerical evidence in the case of a purely harmonic solution.

Different choices of functions for rational approximation were discussed, with a specific treatment for functions in the mass term, that carry a factor $s^2$.
We observed that filtering the unstable poles of the rational approximation is required for obtaining a stable system of differential equations.
We also observed that the Extended AAA method with a degree two polynomial part improved the quality of the approximation.
An additional, but less pronounced, improvement was obtained by dividing the nonlinear term by $s^2$ when a nonlinear function appears in the mass matrix.
The use of Filtered AAA together with Extended AAA gave the most accurate approximation.

The size of the linearization is a multiple of the size of the original system.
For the numerical examples the increase in size was below $10$.
We did not show timings, but, for an implicit time stepper, the cost is dominated by a linear solve with the real valued matrix $R(\sigma)$ for some real $\sigma$,
and matrix vector products with the coefficient matrices.

\bibliographystyle{plain}
\bibliography{stringlong,library}

\end{document}